\documentclass[11pt,english]{article}
\usepackage[T1]{fontenc}
\usepackage[latin2]{inputenc}
\usepackage{a4wide}
\usepackage{amsmath}
\usepackage{color}
\usepackage{graphicx}
\usepackage{amssymb}

\makeatletter

\newcommand{\noun}[1]{\textsc{#1}}

 \newcommand{\lyxaddress}[1]{
   \par {\raggedright #1
   \vspace{1.4em}
   \noindent\par}
 }

\makeatother
\begin{document}

\title{Differential gorms, differential worms}

\author{Denis Kochan, Pavol \v Severa}

\maketitle
\begin{abstract}
We study {}``higher-dimensional'' generalizations of differential
forms. Just as differential forms can be defined as the universal
commutative differential algebra containing $C^{\infty}(M)$, we
can define differential gorms as the universal commutative
bidifferential algebra containing $C^{\infty}(M)$. From a more
conceptual point of view, differential forms are functions on the
superspace of maps $\mathbb{R}^{0|1}\rightarrow M$ and the action
of $\mathit{Diff}(\mathbb{R}^{0|1})$ on forms is equivalent to
deRham differential and to degrees of forms. Gorms are functions
on the superspace of maps $\mathbb{R}^{0|2}\rightarrow M$ and we
study the action of $\mathit{Diff}(\mathbb{R}^{0|2})$ on gorms; it
contains more than just degrees and differentials. By replacing 2
with arbitrary $n$, we get differential worms.

We also study a generalization of homological algebra that one
obtains by replacing $\mathit{Diff}(\mathbb{R}^{0|1})$ with
$\mathit{Diff}(\mathbb{R}^{0|n})$ for $n\ge 2$, and the closely
related question of forms (and gorms and worms) on some
generalized spaces (contravariant functors and stacks) and of
approximations of such {}``spaces'' in terms of worms.

Clearly, this is not a gormless paper.
\end{abstract}
\begin{figure}[h]
\begin{center}\includegraphics[%
  width=6.5cm,
  keepaspectratio]{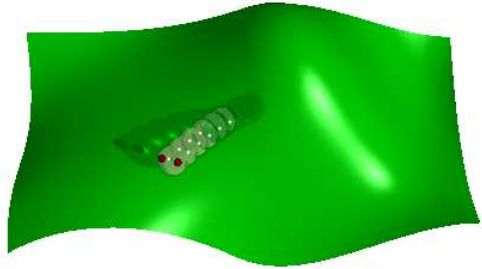}\end{center}

\caption{\label{fig:worm}A worm on a manifold}
\end{figure}

\tableofcontents{}

\section{Introduction}

There is a well known idea of regarding differential forms on a
manifold $M$ as functions on a supermanifold, namely on the odd
tangent bundle $\Pi TM$. The de Rham differential then becomes a
vector field on $\Pi TM$.

This paper is based on the following remarkable fact: one can
describe the supermanifold $\Pi TM$ as the superspace of all maps
$\mathbb{R}^{0|1}\rightarrow M$ (i.e.~as the superspace of all odd
curves in $M$). More importantly, the action of the group
$\mathit{Diff}(\mathbb{R}^{0|1})$ then gives rise to de Rham
differential and to degrees of differential forms. 
This point of view can be found (more or less explicitly)
at many places in physics literature, but we took it explicitly
from \cite{kon}.

We use this idea for two purposes. Firstly, it is natural to make
a generalization and to study the superspaces of all maps $\mathbb{R}^{0|n}\rightarrow M$
for arbitrary $n$, and the action of $\mathit{Diff}(\mathbb{R}^{0|n})$
on these superspaces. To explain the title, the functions on these
map spaces will be called \emph{differential gorms} in the case of
$n=2$, or generally \emph{differential worms} for arbitrary $n$.
Secondly, the idea is straightforwardly applied to some generalizations
of manifolds, namely to contravariant functors from the category of
manifolds, or more generally to stacks. A contravariant functor $F$
from the category of manifolds is usually understood as a {}``generalized
space'', such that $F(M)$ is the set of maps from $M$ to that space.
Differential forms on $F$ should thus be functions on $F(\mathbb{R}^{0|1})$
(and differential gorms functions on $F(\mathbb{R}^{0|2})$, etc.).
As an example, we get a very simple interpretation of equivariant
de Rham theory. The problem of worms of generalized spaces is closely related with a 
generalization of homological algebra, where $\mathit{Diff}(\mathbb{R}^{0|1})$ is replaced
by $\mathit{Diff}(\mathbb{R}^{0|n})$ for arbitrary $n$.

Here is the plan of the paper: In Section 2 we shall describe
differential gorms without use of supermanifolds, as the universal
commutative bidifferential algebra containing the algebra
$C^{\infty}(M)$. This point of view is not completely
satisfactory, as it doesn't reveal the action of the supergroup
$\mathit{Diff}(\mathbb{R}^{0|2})$ (only the sub-supergroup of
affine transformations can be seen). Then, as an appetiser, in
Section 3 we identify differential forms with functions on the
space of all odd curves and derive the basic properties of
differential forms from this fact. In Section 4 we start to
investigate differential gorms (and worms) as functions on the
superspace of all maps $\mathbb{R}^{0|2}\rightarrow M$. In Section
5 we decompose gorms as a representation of
$\mathit{Diff}(\mathbb{R}^{0|2})$ and in Section 6 we prove a
theorem connecting the Euler characteristic with the integral of
any closed integrable gorm. In the final, and possibly the most
interesting Section 7, we look at generalized spaces -
contravariant functors and stacks.

Finally we should add that the ideas used here are very simple; much
of the length of the paper is caused by our attempt to write explicit
coordinate expressions.

\subsubsection*{Remarks on notation}
Generally we denote even coordinates on supermanifolds by latin letters and odd coordinates
by greek letters. To avoid confusion with differential forms, we denote Berezin integral
with respect to (say) $x$ and $\xi$ as $\int f(x,\xi) \overline{dxd\xi}$.

\subsubsection*{Acknowledgement}
We would like to thank to Mari\'an Fecko for many discussions,  to
Peter Pre\v snajder for useful comments on an early version of
this paper, and to James Stasheff for suggesting many improvements
after the first version was posted in the arXives.

\section{Differential gorms as a universal bidifferential algebra}

One can describe the algebra of differential forms $\Omega(M)$ as
the universal graded-commutative differential graded algebra containing
the algebra $C^{\infty}(M)$. That is, there is an algebra homomorphism
$C^{\infty}(M)\rightarrow\Omega^{0}(M)$, and if $\mathcal{A}$ is
any graded-commutative differential graded algebra with a homomorphism
$C^{\infty}(M)\rightarrow\mathcal{A}^{0}$, there is a unique homomorphism
of differential graded algebras $\Omega(M)\rightarrow\mathcal{A}$
making the triangle commutative.

In the same spirit, we can define the algebra of \emph{differential
gorms} $\Omega_{[2]}(M)$ as the universal graded-commutative bidifferential
algebra containing $C^{\infty}(M)$ (by a bidifferential algebra we
mean a bicomplex with a compatible structure of algebra; in particular,
it is $\mathbb{Z}^{2}$-graded). That is, there is an algebra homomorphism
$C^{\infty}(M)\rightarrow\Omega_{[2]}^{0,0}(M)$, and if $\mathcal{A}$
is any graded-commutative bidifferential algebra with a homomorphism
$C^{\infty}(M)\rightarrow\mathcal{A}^{0,0}$, there is a unique homomorphism
of bidifferential algebras $\Omega_{[2]}(M)\rightarrow\mathcal{A}$
making the triangle commutative. Just as in the case of $\Omega(M)$
it turns out that $\Omega_{[2]}^{0,0}(M)\cong C^{\infty}(M)$.

To make this abstract definition down-to-earth, let us choose
local coordinates $x^{i}$ on $M$. The algebra $\Omega_{[2]}(M)$ is
freely generated by the algebra of functions, and by the elements
$d_{1}x^{i}$, $d_{2}x^{i}$ and $d_{1}d_{2}x^{i}$ (of bidegrees
$(1,0)$, $(0,1)$ and $(1,1)$ respectively), where $d_{a}$
($a=1,2$) are the two differentials in $\Omega_{[2]}(M)$. If $f$
is a function then $d_{a}f=\frac{\partial f}{\partial
x^{i}}d_{a}x^{i}$ and therefore
$d_{1}d_{2}f=\frac{\partial^{2}f}{\partial x^{i}\partial
x^{j}}d_{1}x^{i}d_{2}x^{j}+\frac{\partial f}{\partial
x^{i}}d_{1}d_{2}x^{i}$. If $\tilde{x}^{i}$ is another system of
local coordinates, we get from here\[
d_{a}\tilde{x}^{i}=\frac{\partial\tilde{x}^{i}}{\partial
x^{j}}d_{a}x^{j}\]
\[
d_{1}d_{2}\tilde{x}^{i}=\frac{\partial^{2}\tilde{x}^{i}}{\partial x^{j}\partial x^{k}}d_{1}x^{j}d_{2}x^{k}+\frac{\partial\tilde{x}^{i}}{\partial x^{j}}d_{1}d_{2}x^{j}.\]
 Differential gorms are not tensor fields, since their transformation
law involves second derivatives; they belong to 2nd order geometry.

Finally, for arbitrary $n$ we can define the algebra of \emph{differential
worms of level n,} $\Omega_{[n]}(M)$, as the universal graded-commutative
$n$-differential algebra containing $C^{\infty}(M)$. One can easily
compute transformation laws for worms of any level; they contain $n$-th
derivatives at most.

On $\Omega_{[n]}(M)$ we have an obvious action of the semigroup $\mathit{Mat}(n)$:
it leaves $C^{\infty}(M)$ intact, and linearly transforms the $n$
differentials $d_{a}$. It turns out that the differentials and the
action of $\mathit{Mat}(n)$ are just the tip of an iceberg: on $\Omega_{[n]}(M)$
we have an action of the supersemigroup of all maps $\mathbb{R}^{0|n}\rightarrow\mathbb{R}^{0|n}$
($\mathit{Mat}(n)$ corresponds to linear transformations of $\mathbb{R}^{0|n}$
and the differentials to translations). To see this will require a
new point of view on differential worms.

\section{Appetizer: differential forms}

\subsection{$\Pi TM$ as the space of odd curves}

We denote by $\Pi TM$ the supermanifold $M^{\mathbb{R}^{0|1}}$ of
all maps $\mathbb{R}^{0|1}\rightarrow M$. It is characterized by
the following property: for any supermanifold $Y$, a map $Y\rightarrow\Pi TM$
is the same as a map $\mathbb{R}^{0|1}\times Y\rightarrow M$ (in
other words, the functor $\Pi T$ is the right adjoint of the functor
$Y\mapsto\mathbb{R}^{0|1}\times Y$).

It is easy to understand $\Pi TM$ in local coordinates. If $\theta$
is the coordinate on $\mathbb{R}^{0|1}$ and $x^{i}$ are local coordinates
on $M$, a map $\mathbb{R}^{0|1}\rightarrow M$ parametrized by $Y$,
i.e. a map $\mathbb{R}^{0|1}\times Y\rightarrow M$, is given by functions
\[
x^{i}(\theta,\eta)=x^{i}(\eta)+\theta\xi^{i}(\eta),\]
where $\eta$ denotes local coordinates on $Y$ (we just used Taylor
expansion in $\theta$). Such a map is therefore the same as a map
from $Y$ to a supermanifold with even coordinates $x^{i}$ and odd
coordinates $\xi^{i}$. To find $\Pi TM$ globally, suppose $\tilde{x}^{i}$
is another system of local coordinates on $M$; then \[
\tilde{x}^{i}(\theta)=\tilde{x}^{i}(x(\theta))=\tilde{x}^{i}(x+\theta\xi)=\tilde{x}^{i}(x)+\theta\frac{\partial\tilde{x}^{i}}{\partial x^{j}}\xi^{j},\]
i.e.~the transition functions on $\Pi TM$ are\[
\tilde{x}^{i}=\tilde{x}^{i}(x),\;\tilde{\xi}^{i}=\frac{\partial\tilde{x}^{i}}{\partial x^{j}}\xi^{j}.\]
Now we can see that we can identify functions on $\Pi TM$ with differential
forms on $M$, by identifying $\xi^{i}$ with $dx^{i}$ (we also see
that our notation is correct, i.e.~that the supermanifold is indeed
the odd tangent bundle of $M$).

\subsection{De Rham differential and action of $\mathit{Diff}(\mathbb{R}^{0|1})$}

Since $\Pi TM$ is the supermanifold of all maps $\mathbb{R}^{0|1}\rightarrow M$,
we have an action of the supergroup $\mathit{Diff}(\mathbb{R}^{0|1})$
on it. Let us first describe $\mathit{Diff}(\mathbb{R}^{0|1})$. It
is an open sub-supergroup of the supersemigroup $\Pi T\mathbb{R}^{0|1}$
of all maps $\mathbb{R}^{0|1}\rightarrow\mathbb{R}^{0|1}$. A (possibly
parametrized) map $\mathbb{R}^{0|1}\rightarrow\mathbb{R}^{0|1}$ is
of the form $\theta'=a\theta+\beta$, hence the supersemigroup is
diffeomorphic to $\mathbb{R}^{1|1}$, with one even coordinate $a$
and one odd coordinate $\beta$; $\mathit{Diff}(\mathbb{R}^{0|1})$
is given by $a\neq0$.

The right action of $\Pi T\mathbb{R}^{0|1}$ (hence also of $\mathit{Diff}(\mathbb{R}^{0|1})$)
on $\Pi TM$ is given by $x'+\theta\xi'=x+(a\theta+\beta)\xi$, i.e.\[
x'=x+\beta\xi,\;\xi'=a\xi.\]
 The vector fields generating this action are\[
E=\xi^{i}\frac{\partial}{\partial\xi^{i}},\; d=\xi^{i}\frac{\partial}{\partial x^{i}}.\]

When $E$ acts on a function on $\Pi TM$, i.e.~on a differential
form on $M$, it multiplies it by its degree; on the other hand, $d$
acts as the de Rham differential. The canonical structure of a complex
on $\Omega(M)$ is therefore equivalent to the action of $\mathit{Diff}(\mathbb{R}^{0|1})$.
The fact that it is a complex follows from the commutation relations\[
[E,d]=d,\;[d,d]=0\]
 in the Lie algebra of $\mathit{Diff}(\mathbb{R}^{0|1})$.

\subsection{Cartan formula and its generalizations\label{sub:De-Rham-differential}}

Now we will try to understand Cartan's formula $\mathcal{L}_{v}=d\, i_{v}+i_{v}\, d$
from this point of view. It will be convenient to treat a general
map space $Y^{X}$ of all maps $X\rightarrow Y$ between two (super)manifolds;
we shall however pretend that $Y^{X}$ is finite-dimensional (since
it is so for $X=\mathbb{R}^{0|n}$) to avoid analytical problems. As we shall see, 
Cartan's formula comes from the action of the group $\mathit{Diff}(X)\ltimes\left(\mathit{Diff}(Y)\right)^{X}$ on $Y^X$ (in the case of $X=\mathbb{R}^{0|1}$).

Let $\mathit{e}v:X\times Y^{X}\rightarrow Y$ be the evaluation map
(adjoint to the identity map $Y^{X}\rightarrow Y^{X}$); in the case
of $X=\mathbb{R}^{0|1}$ and $Y=M$ (so that $Y^{X}=\Pi TM$), $\mathit{ev}:(\theta,x,\xi)\mapsto x+\theta\xi.$
We can naturally identify vector fields on $Y^{X}$ with sections
of the vector bundle $\mathit{ev}^{*}TY$,%
\footnote{to see this, realize that a vector field on $Y^{X}$ is an infinitesimal
deformation of the identity map $Y^{X}\rightarrow Y^{X}$, i.e.~an
infinitesimal deformation of $\mathit{e}v:X\times Y^{X}\rightarrow Y$,
i.e~a section of $\mathit{ev}^{*}TY$%
} and therefore we can multiply them with arbitrary functions on $X\times Y^{X}$
(not just on $Y^{X}$). If $f$ is a function on $X$ (and therefore
also on $X\times Y^{X}$) and $w$ a vector field on $Y^{X}$, we
shall denote their product as $f\cdot w$.

If $u$ is a vector field on $X$ and $v$ a vector field on $Y$,
we shall denote their natural lifts to $Y^{X}$ as $u^{\flat}$ and
$v^{\sharp}$; these are the vector fields generating the \emph{left}
actions of $\mathit{Diff}(X)$ and $\mathit{Diff}(Y)$ on $Y^{X}$.
The following formulas express the fact that on $Y^{X}$ we have a
(left) action of the semidirect product $\mathit{Diff}(X)\ltimes\left(\mathit{Diff}(Y)\right)^{X}$:\[
[u_{1}^{\flat},u_{2}^{\flat}]=[u_{1},u_{2}]^{\flat}\]
\[
[u^{\flat},f\cdot v^{\sharp}]=(uf)\cdot v^{\sharp}\]
\[
[f_{1}\cdot v_{1}^{\sharp},f_{2}\cdot v_{2}^{\sharp}]=(-1)^{|f_{2}||v_{1}|}f_{1}f_{2}\cdot[v_{1},v_{2}]^{\sharp};\]
we also have\[
f\cdot u^{\flat}=(fu)^{\flat}.\]

In the case of $X=\mathbb{R}^{0|1}$ and $Y=M$ we have $(\partial_{\theta})^{\flat}=-d$,
$(\theta\partial_{\theta})^{\flat}=-E$ (the minus signs are here
since $d$ and $E$ generate the \emph{right} action), $\theta\cdot\partial_{x^{i}}=-\partial_{\xi^{i}}$
(and therefore $\theta\cdot\partial_{\xi^{i}}=0$). Moreover, for
any vector field $v$ on $M$ we have $v^{\sharp}=\mathcal{L}_{v}$
and $\theta\cdot v^{\sharp}=-i_{v}$. The equation $[(\partial_{\theta})^{\flat},\theta\cdot v^{\sharp}]=v^{\sharp}$
is Cartan's $d\, i_{v}+i_{v}\, d=\mathcal{L}_{v}$.

\subsection{Vector fields on $\Pi TM$\label{sub:Vector-fields-on}}

As we have seen, the space of vector fields on $Y^{X}$ is a module
over $C^{\infty}(X)$, and it is also a representation of $\mathit{Diff}(X)$.
If $\phi\in\mathit{Diff}(X)$, $f\in C^{\infty}(X)$ and $w$ is a
vector field on $Y^{X}$ then clearly\[
\phi\cdot(f\cdot w)=(f\circ\phi^{-1})\cdot(\phi\cdot w);\]
 in other words, the space of vector fields on $Y^{X}$ is a module
over the crossed product of $C^{\infty}(X)$ with $\mathit{Diff}(X)$.
Infinitesimally, if $u$ is a vector field on $X$,\[
[u^{\flat},f\cdot w]=(uf)\cdot w+(-1)^{|f||u|}f\cdot[u^{\flat},w].\]

Now let us return to the case of $X=\mathbb{R}^{0|1}$, $Y=M$, $Y^{X}=\Pi TM$.
Vector fields on $Y^{X}$ are then derivations of the algebra $\Omega(M)$.
Notice that $\theta\cdot(\theta\cdot w)=0$ (since $\theta^{2}=0$)
and that $[(\partial_{\theta})^{\flat},[(\partial_{\theta})^{\flat},w]]=0$
(since $[(\partial_{\theta})^{\flat},(\partial_{\theta})^{\flat}]=0$),
i.e.~both $\theta$ and $(\partial_{\theta})^{\flat}$ act as differentials;
$(\partial_{\theta})^{\flat}$ increases degree by 1, while $\theta$
decreases it by 1. Finally, \[
w=[(\partial_{\theta})^{\flat},\theta\cdot w]+\theta\cdot[(\partial_{\theta})^{\flat},w],\]
i.e.~any $w$ can be \emph{uniquely} decomposed as $w=w_{1}+w_{2}$
(by $w_{1}=[(\partial_{\theta})^{\flat},\theta\cdot w]$, $w_{2}=\theta\cdot[(\partial_{\theta})^{\flat},w]$)
so that $[(\partial_{\theta})^{\flat},w_{1}]=0$ and $\theta\cdot w_{2}=0$.

In coordinates, a vector field $w$ on $\Pi TM$, of degree $p$ and
such that $\theta\cdot w=0$, is of the form \[
w=A_{i_{1}i_{2}\dots i_{p+1}}^{k}\xi^{i_{1}}\xi^{i_{2}}\dots\xi^{i_{p+1}}\partial_{\xi^{k}},\]
i.e.~it is a section of $TM\otimes\bigwedge^{p+1}T^{*}M$.

Let us summarize what we have found. The graded Lie algebra $\mathit{Der}(\Omega(M))$
decomposes to a direct sum of graded vector spaces\[
\mathit{Der}(\Omega(M))=\ker(\theta)\oplus\ker(\partial_{\theta})\]
(they turn out to be subalgebras). The two subspaces are naturally
isomorphic; the two mutually inverse isomorphisms are the action of
$\theta$ ($\ker(\partial_{\theta})\rightarrow\ker(\theta)$) and
the action of $\partial_{\theta}$ ($\ker(\theta)\rightarrow\ker(\partial_{\theta})$).
Moreover, they are both isomorphic with the space of vector fields
with values in differential forms. The Lie bracket on $\ker(\partial_{\theta})$
(the derivations of $\Omega(M)$ commuting with the differential)
is the Fr\"olicher-Nijenhuis bracket.

\subsection{Integration and Stokes formula}

The fact that $\Pi TM$ is a map space doesn't seem to shed much
light on integration of differential forms, i.e.~of functions on
$\Pi TM$. For this reason we just repeat the simple standard
facts: the volume measure $\overline{dxd\xi}$ on $\Pi TM$ is
independent of the choice of coordinate system $x^i$ (up to choice
of orientation), and it is also invariant with respect to the
vector field $d$. For that reason $\int d\alpha=0$ for any form
$\alpha$ with compact support. To get Stokes theorem, let
$\chi_{\Omega}$ be the characteristic function of a compact domain
$\Omega$; then $0=\int
d(\chi_{\Omega}\alpha)=\int(d\chi_{\Omega})\alpha+\int\chi_{\Omega}d\alpha=-\int_{\partial\Omega}\alpha+\int_{\Omega}d\alpha$.

It is certainly an interesting thing that on the space of all maps
$\mathbb{R}^{0|1}\rightarrow M$ there is a natural volume measure
(e.g.~from the point of view of quantum field theory). As we will
see, the situation gets even better for $\mathbb{R}^{0|n}$ with $n\ge2$;
the measure will be $\mathit{Diff}(\mathbb{R}^{0|n})$-invariant (here
it was $d$-invariant, but not $E$-invariant).

\section{Differential gorms as functions on the space of odd surfaces }

Everything we'll be doing here will be fairly analogous to the previous
section, so we can be brief. Let $(\Pi T)^{2}M$ denote the supermanifold
of all maps $\mathbb{R}^{0|2}\rightarrow M$. Let $\theta^{1}$, $\theta^{2}$
be the coordinates on $\mathbb{R}^{0|2}$ and $x^{i}$ be local coordinates
on $M$. A (parametrized) map $\mathbb{R}^{0|2}\rightarrow M$ expanded
to Taylor series in $\theta$'s looks as \[
x^{i}(\theta^{1},\theta^{2})=x^{i}+\theta^{1}\xi_{1}^{i}+\theta^{2}\xi_{2}^{i}+\theta^{2}\theta^{1}y^{i},\]
 $(\Pi T)^{2}M$ has therefore local coordinates $x^{i}$, $y^{i}$
(even coordinates) and $\xi_{1}^{i}$, $\xi_{2}^{i}$ (odd coordinates).
If $\tilde{x}^{i}$ is another system of local coordinates on $M$
then (expanding to Taylor series)\[
\tilde{x}^{i}(x(\theta^{1},\theta^{2}))=\tilde{x}^{i}(x+\theta^{1}\xi_{1}+\theta^{2}\xi_{2}+\theta^{2}\theta^{1}y)=\]
\[
=\tilde{x}^{i}(x)+\theta^{1}\frac{\partial\tilde{x}^{i}}{\partial x^{j}}\xi_{1}^{j}+\theta^{2}\frac{\partial\tilde{x}^{i}}{\partial x^{j}}\xi_{2}^{j}+\theta^{2}\theta^{1}\left(\frac{\partial\tilde{x}^{i}}{\partial x^{j}}y^{j}+\frac{\partial^{2}\tilde{x}^{i}}{\partial x^{j}\partial x^{k}}\xi_{1}^{j}\xi_{2}^{k}\right),\]
 i.e.~the transition functions on $(\Pi T)^{2}M$ are \[
\tilde{x}^{i}=\tilde{x}^{i}(x),\;\tilde{\xi}_{a}^{i}=\frac{\partial\tilde{x}^{i}}{\partial x^{j}}\xi_{a}^{j},\;\tilde{y}^{i}=\frac{\partial\tilde{x}^{i}}{\partial x^{j}}y^{j}+\frac{\partial^{2}\tilde{x}^{i}}{\partial x^{j}\partial x^{k}}\xi_{1}^{j}\xi_{2}^{k}.\]
 From this we can see that we can identify differential gorms on $M$
(as defined in section 2) with functions on $(\Pi T)^{2}M$ polynomial
in $y$'s, by identifying $d_{a}x^{i}$ with $\xi_{a}^{i}$ and $d_{1}d_{2}x^{i}$
with $y^{i}$. We arrived to this identification by a computation,
but there is a simpler reason using differentials (see below). General
functions on $(\Pi T)^{2}M$ (not necessarily polynomial in $y$'s)
will be called \emph{pseudodifferential gorms} (these are things like
$e^{-(d_{1}d_{2}x)^{2}}d_{1}x\, d_{2}x$).

Let us notice that the body of the supermanifold $(\Pi T)^{2}M$ is
naturally isomorphic to $TM$; indeed, we get the body by setting
all the odd coordinates $\xi_{a}^{i}$ to zero, and then the transition
functions for $x$'s and $y$'s become those of $TM$.

Finally, let us describe $(\Pi T)^{n}M=M^{\mathbb{R}^{0|n}}$ for
higher $n$'s and directly identify its functions with differential
worms of level $n$ on $M$. If $x^{i}$'s are local coordinates on
$M$ then the coordinates on $(\Pi T)^{n}M$ are $x^{i}$, $d_{a}x^{i}$
$(1\leq a\leq n)$, $d_{a}d_{b}x^{i}$ $(1\leq a<b\leq n)$, \ldots{}
, $d_{1}d_{2}\dots d_{n}x^{i}$ (i.e.~apply the differentials $d_{a}$
$(1\leq a\leq n)$ to $x^{i}$'s in all possible ways). These coordinates
are identified with Taylor coefficients of a map $\mathbb{R}^{0|n}\rightarrow M$
by\[
x^{i}(\theta^{1},\theta^{2},\dots,\theta^{n})=e^{\theta^{a}d_{a}}x^{i}.\]
From this expression it is clear that $-\left(\partial_{\theta^{a}}\right)^{\flat}$
is equal to $d_{a}$. For example, in the case of $n=2$ we have $-(\partial_{\theta^{1}})^{\flat}=\xi_{1}^{i}\partial_{x^{i}}+y^{i}\partial_{\xi_{2}^{i}}$,
$-(\partial_{\theta^{2}})^{\flat}=\xi_{2}^{i}\partial_{x^{i}}-y^{i}\partial_{\xi_{1}^{i}}$,
i.e.~$-(\partial_{\theta^{a}})^{\flat}=\xi_{a}^{i}\partial_{x^{i}}+\epsilon_{ab}y^{i}\partial_{\xi_{b}^{i}}$.

\subsection{More than differentials: action of $\mathit{Diff}(\mathbb{R}^{0|2})$}

We have already found the vector fields $-\left(\partial_{\theta^{a}}\right)^{\flat}$:
they generate the action of the group of translations of $\mathbb{R}^{0|2}$
on $(\Pi T)^{2}M$, and they are equal to the two differentials on
differential gorms. We can easily find $u^{\flat}$ for any vector
field $u$ on $\mathbb{R}^{0|2}$: either we do it directly, regarding
$u$ as an infinitesimal transformation on $\mathbb{R}^{0|2}$ and
finding the corresponding infinitesimal transformation of $(\Pi T)^{2}M$
from the formula\begin{equation}
x'+\theta^{a}\xi'_{a}+\theta^{2}\theta^{1}y'=x+\theta'^{a}\xi_{a}+\theta'^{2}\theta'^{1}y,\label{eq:posobenie}\end{equation}
or use the known expression for $\left(\partial_{\theta^{a}}\right)^{\flat}$
and the identities \[
\theta^{a}\cdot\partial_{x^{i}}=-\partial_{\xi_{a}^{i}},\;\;\theta^{2}\theta^{1}\cdot\partial_{x^{i}}=\partial_{y^{i}}.\]
The result is\[
d_{a}=-(\partial_{\theta^{a}})^{\flat}=\xi_{a}^{i}\partial_{x^{i}}+\epsilon_{ab}y^{i}\partial_{\xi_{b}^{i}}\]
\[
E_{a}^{b}=-\left(\theta^{b}\partial_{\theta^{a}}\right)^{\flat}=\xi_{a}^{i}\partial_{\xi_{b}^{i}}+\delta_{b}^{a}y^{i}\partial_{y^{i}}\]
\[
R_{a}=-(\theta^{2}\theta^{1}\partial_{\theta^{a}})^{\flat}=\xi_{a}^{i}\partial_{y^{i}}.\]
These vector fields give us a right action of the Lie algebra $\mathfrak{diff}(\mathbb{R}^{0|2})$
on $(\Pi T)^{2}M$ (that is the reason for the minus signs: $u^{\flat}$
give the left action), i.e.~a left action on the algebra of differential
gorms. $E_{1}^{1}$ and $E_{2}^{2}$ are the two degrees on differential
gorms, $E_{a}^{b}$ generate the action of $\mathfrak{gl}(2)$ on
gorms, but $R_{a}$'s are something new. Similar formulas can be easily
found for worms of arbitrary level.

If we want to know the right action of the supersemigroup $(\Pi T)^{2}\mathbb{R}^{0|2}=\left(\mathbb{R}^{0|2}\right)^{\mathbb{R}^{0|2}}$on
$(\Pi T)^{2}M$ (not just the infinitesimal action we have just derived),
we can easily find it from (\ref{eq:posobenie}). A (parametrized)
map $\mathbb{R}^{0|2}\rightarrow\mathbb{R}^{0|2}$ is of the form
\[
\theta'^{1}=\beta^{1}+a_{1}^{1}\theta^{1}+a_{2}^{1}\theta^{2}+\gamma^{1}\theta^{2}\theta^{1}\]
\[
\theta'^{2}=\beta^{2}+a_{1}^{2}\theta^{1}+a_{2}^{2}\theta^{2}+\gamma^{2}\theta^{2}\theta^{1},\]
 i.e.~$(\Pi T)^{2}\mathbb{R}^{0|2}$ is diffeomorphic to $\mathbb{R}^{4|4}$
(with even coordinates $a_{b}^{a}$ and odd coordinates $\beta^{a}$,
$\gamma^{a}$); we also see that its body is the semigroup $\mathit{Mat}(2)$
of $2\times2$-matrices. The action of $\mathit{Mat}(2)$ (i.e.~when
we set $\beta$'s and $\gamma$'s to zero) is given by\[
x'=x,\;\xi'_{a}=a_{a}^{b}\xi_{b},\; y'=\det(A)y,\]
 where $A$ is the matrix $a_{b}^{a}$. If it ever becomes useful,
the full result of (\ref{eq:posobenie}) is

\[
x'=x+\beta^{a}\xi_{a}+\beta^{2}\beta^{1}y,\;\xi'_{a}=a_{a}^{b}\xi_{b}+\epsilon_{bc}\beta^{b}a_{a}^{c}y,\;
y'=(\det(A)+\epsilon_{bc}\beta^{b}\gamma^{c})y+\gamma^{a}\xi_{a}.\]

\subsection{Cartan gormulas}

Analogues of the Cartan formula $\mathcal{L}_{v}=d\, i_{v}+i_{v}\, d$
(i.e.~$v^{\sharp}=[(\partial_{\theta})^{\flat},\theta\cdot v^{\sharp}]$)
for worms of arbitrary level were already found in section \ref{sub:De-Rham-differential}.
Here we mention just one special case: since $\theta^{2}\theta^{1}\cdot\partial_{x^{i}}=\partial_{y^{i}}$,
we have $\theta^{2}\theta^{1}\cdot v^{\sharp}=v^{i}\partial_{y^{i}}$
for any vector field $v=v^{i}\partial_{x^{i}}$ on $M$; let us denote
$\theta^{2}\theta^{1}\cdot v^{\sharp}$ as $i_{v}$. Then $[d_{1},[d_{2},i_{v}]]=\mathcal{L}_{v}$.

\subsection{Cohomology of gorms}

The cohomology of differential gorms with respect to any $d_{a}$
is naturally isomorphic to de Rham cohomology of $M$. To see this
we just write $(\Pi T)^{2}M$ as $\Pi T(\Pi TM)$, i.e.~the cohomology
of differential gorms is de Rham cohomology of $\Pi TM$, and $\Pi TM$
can be contracted to $M$.

We can also compute cohomology with respect to say $R_{1}$, since
$(R_{1})^{2}=0$. This cohomology is isomorphic to $\Omega(M)$. In
fact, the projection $\Omega_{[2]}(M)\rightarrow\Omega(M)$ given
by $d_{1}\mapsto d$, $d_{2}\mapsto0$, is a quasiisomorphism, when
$\Omega(M)$ is taken with zero differential and $\Omega_{[2]}(M)$
with differential $R_{1}$. The reason is simple: $[R_{1},d_{2}]=-E_{1}^{1}$,
i.e.~$E_{1}^{1}$ acts trivially on cohomology, i.e.~the semigroup
action $x\mapsto x$, $d_{1}x\mapsto\lambda d_{1}x$, $d_{2}x\mapsto d_{2}x$,
$d_{1}d_{2}x\mapsto\lambda d_{1}d_{2}x$ is trivial on cohomology,
and setting $\lambda=0$ gives us the result. Notice that we could
use the same argument to show that the cohomology with respect to
$d_{2}$ is just the cohomology of $\Omega(M)$ with respect to $d$.
In geometric terms, we have an embedding $\Pi TM\subset(\Pi TM)^{2}$
coming from the projection $\mathbb{R}^{0|2}\rightarrow\mathbb{R}^{0|1}$,
$(\theta^{1},\theta^{2})\mapsto\theta^{2}$; this projection can be
obtained by the semigroup action $(\theta^{1},\theta^{2})\mapsto(\lambda\theta^{1},\theta^{2})$,
setting $\lambda=0$. The action is generated by $\theta^{1}\partial_{\theta^{1}}$
and $\theta^{1}\partial_{\theta^{1}}=[\theta^{2}\theta^{1}\partial_{\theta^{1}},\partial_{\theta^{2}}]$.

\subsection{Integration of gorms\label{sub:Integration-of-gorms}}

Like differential forms, gorms (and worms) can be integrated. However
we need pseudodifferential gorms (depending non-polynomially on $d_{1}d_{2}x$'s)
for the integral to be finite. In coordinates, the integral is just
the Berezin integral with volume measure $\overline{dx^{1}\, d\xi_{1}^{1}\, d\xi_{2}^{1}\, dy^{1}\dots dx^{m}\, d\xi_{1}^{m}\, d\xi_{2}^{m}\, dy^{m}}$,
where $m$ is the dimension of $M$. In other words, to integrate
a gorm, expand it in $\xi$'s, take the coefficient in front of $\xi_{1}^{1}\xi_{2}^{1}\xi_{1}^{2}\xi_{2}^{2}\dots\xi_{1}^{m}\xi_{2}^{m}$,
and integrate it over $y$'s and $x$'s. For example, if $M=\mathbb{R}$,
\[
\int e^{-x^{2}-(d_{1}d_{2}x)^{2}}d_{1}x\, d_{2}x=\pi.\]

The integral is independent of the choice of coordinates, i.e.~$\overline{dx\, d\xi_{1}\, d\xi_{2}\, dy}$
is $\mathit{Diff}(M)$-invariant, and in fact it is also $\mathit{Diff}(\mathbb{R}^{0|2})$-invariant.
It means that for any integrable gorm $\alpha$ and any vector field
$u$ on $\mathbb{R}^{0|2}$ we have $\int(u^{\flat}\alpha)=0$ --
a form of Stokes theorem.

A similar claim is true for worms with any $n\geq2$; the
coordinate Berezin integral on $(\Pi T)^{n}M$ is both
$\mathit{Diff}(\mathbb{R}^{0|n})$ and
$\mathit{Diff}(M)$-invariant.


\section{Differential gorms as a representation of $\mathit{Diff}(\mathbb{R}^{0|2})$\label{sec:representation}}

Our aim in this section is to decompose $\Omega_{[2]}(M)$ to
indecomposable representations of
$\mathit{Diff}(\mathbb{R}^{0|2})$. It should be compared with the
de Rham cohomology of $\Omega(M)$; the latter is connected with
the problem "solve the equation $d\alpha=\beta$  in
$\Omega(M)$", while the decomposition of $\Omega_{[2]}(M)$
describes solutions of all linear equations in $\Omega_{[2]}(M)$
that use the operators $d_a$, $E^b_a$ and $R_a$ (i.e.~the action of $\mathit{Diff}(\mathbb{R}^{0|2})$) on $\Omega_{[2]}(M)$.

The representation theory of $\mathit{Diff}(\mathbb{R}^{0|n})$ for
$n\ge 3$ was shown to be wild by N.~Shomron \cite{sho} (i.e.~it is impossible
to classify all finite-dimensional indecomposable representations
of this supergroup). The decomposition of $\Omega_{[n]}$ for $n\ge
3$ remains an open problem for us.

\subsection{Irreducible representations of $\mathit{Mat}(n)$ and of the categories
$\mathit{Vect}$ and $\mathit{Diff}^{op}$}

Let us recall that irreducible representations of the semigroup $\mathit{Mat}(n)$,
or more invariantly, of the semigroup $\mathit{End}(V)$, where $V$
is an $n$-dimensional vector space, are classified by highest weights,
that is by $n$-tuples of integers $l_{1}\geq l_{2}\geq\dots\geq l_{n}\geq0$.
Such an $n$-tuple can be represented by a Young table; for example,
the triple $(5,2,1)$ is represented by

\begin{center}\includegraphics[%
  scale=0.3]{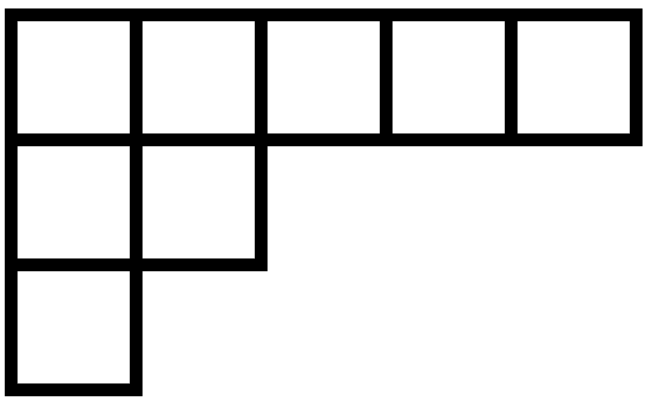}\end{center}

The representation with the highest weight $\vec{l}=(l_{1},l_{2},\dots,l_{n})$
can be found in the decomposition of the representation $V^{\otimes N}$
to irreducibles, where $N=l_{1}+l_{2}+\dots+l_{n}$. Namely, let $W_{\vec{l}}$
be the irreducible representation of the symmetric group $S_{N}$,
corresponding to the Young table $\vec{l}$. $S_{N}$ acts also on
$V^{\otimes N}$, by permutations. The irreducible representation
of $\mathit{End}(V)$ with the highest weight $\vec{l}$ is then \[
V_{\vec{l}}=\mathit{Hom}_{S_{N}}(W_{\vec{l}},V^{\otimes N}).\]
In fact, this formula gives us a representation of the category $\mathit{Vect}$
of finite-dimensional vector spaces, i.e.~a functor $\mathit{Vect}\rightarrow\mathit{Vect}$,
$V\mapsto V_{\vec{l}}$.

For any manifold $M$ we can consider the space of sections of the
vector bundle $T_{\vec{l}}^{*}M$; it is a right representation of
the semigroup of all smooth maps $M\rightarrow M$. This bundle is
non-zero iff the number of rows of $\vec{l}$ is at most $\dim M$,
and the representation is known to be irreducible unless $\vec{l}$
has only one column; in that case, $T_{\vec{l}}^{*}M=\bigwedge^{N}T^{*}M$,
where $N$ is the length of the column, and the space of differential
$N$-forms on $M$ has the invariant subspace of all closed $N$-forms.
Generally, complete reducibility doesn't hold for the representations
of these semigroups; we shall meet many examples soon.

Let is notice that the construction $M\mapsto\Gamma(T_{\vec{l}}^{*}M)$
is a contravariant functor from the category $\mathit{Diff}$ of smooth
manifolds to the category of vector spaces. The functor $\Gamma(T_{\vec{l}}^{*})$
can also be applied to supermanifolds. Let us however notice that
if $X$ is a supermanifold, the action of $S_{N}$ on $(T^{*})^{\otimes N}X$
is modified by the sign rule: the transposition acts by $a\otimes b\mapsto(-1)^{|a||b|}b\otimes a$.
This implies that $T_{\vec{l}}^{*}\mathbb{R}^{0|n}$ is zero iff $\vec{l}$
has more than $n$ \emph{columns}. For example, there is no Riemann
metric on $\mathbb{R}^{0|1}$ (the Young table is $\square\!\square$),
but there are non-zero $k$-forms for any $k$.

$\Omega_{[2]}$ is also a right representation of the category $\mathit{Diff}$
(i.e.~a contravariant functor from $\mathit{Diff}$), and it is also
a representation of the semigroup $\left(\mathbb{R}^{0|2}\right)^{\mathbb{R}^{0|2}}$.
Our aim will be to decompose it to indecomposable parts.

\subsection{Decomposition of $\Omega_{[2]}(M)$}

\subsubsection{Decomposition to $\mathit{Mat}(2)$ irreducibles}

 In this preliminary section we decompose $\Omega_{[2]}(M)$ as a representation of $\mathit{Mat}(2)\subset\left(\mathbb{R}^{0|2}\right)^{\mathbb{R}^{0|2}}$.
It contains only the representations with highest weights $(l_{1},l_{2})$
such that $l_{1}-l_{2}\leq\dim M$; the picture looks something like
this (for $\dim M=4$):\[
\qquad\qquad\qquad\qquad\qquad\qquad\mathit{etc}\dots\]
\begin{equation}
l_{2}\uparrow\begin{array}{ccccccccccc}
5 &  &  &  &  &  & \bullet & \bullet & \bullet & \bullet & \bullet\\
4 &  &  &  &  & \bullet & \bullet & \bullet & \bullet & \bullet\\
3 &  &  &  & \bullet & \bullet & \bullet & \bullet & \bullet\\
2 &  &  & \bullet & \bullet & \bullet & \bullet & \bullet\\
1 &  & \bullet & \bullet & \bullet & \bullet & \bullet\\
0 & \bullet & \bullet & \bullet & \bullet & \bullet\\
 & 0 & 1 & 2 & 3 & 4 & 5 & 6 & 7 & 8 & 9\end{array}\label{eq:kruzky}\end{equation}
\[
l_{1}\rightarrow\]
These representations can be found by looking for gorms with highest
weights. For example, a general gorm of weight (or bidegree) $(2,1)$
is of the form $a_{ij}(x)\, y^{i}\xi_{1}^{j}+b_{ijk}(x)\,\xi_{1}^{i}\xi_{1}^{j}\xi_{2}^{k}$;
the highest weight condition means that it is annulled by $E_{1}^{2}=\xi_{1}^{i}\partial_{\xi_{2}^{i}}$,
in this case it means that $b_{ijk}(x)\,\xi_{1}^{i}\xi_{1}^{j}\xi_{1}^{k}=0$,
i.e.~$b_{ijk}$ becomes 0 after complete skew-symmetrization in $ijk$.

\subsubsection{Generic part: the bundles $\widetilde{T_{\vec{l}}^{*}}M$ (cotangent
tetris)}

Now we would like to decompose $\Omega_{[2]}(M)$ as a representation
of $\left(\mathbb{R}^{0|2}\right)^{\mathbb{R}^{0|2}}$ (for a review
of the needed representation theory see \cite{lei}). Let us first
describe the \emph{right} irreducible representations of $\left(\mathbb{R}^{0|2}\right)^{\mathbb{R}^{0|2}}$;
the left irreducible representations are their duals. For any Young
table $\vec{l}$ with two columns the representation $\Gamma(T_{\vec{l}}^{*}\mathbb{R}^{0|2})$
is irreducible; these are called \emph{generic} irreducibles. The
remaining irreducibles are the spaces of closed differential $k$-forms
on $\mathbb{R}^{0|2}$ for any $k$. The representation theory of
$\left(\mathbb{R}^{0|2}\right)^{\mathbb{R}^{0|2}}$ is quite simple:
the generic irreducibles can't appear in the composition series of
any reducible indecomposable representation.

There is a simple geometrical (and somewhat tautological) way to get
an intertwining map $\Gamma(T_{\vec{l}}^{*}\mathbb{R}^{0|2})^{*}\rightarrow\Omega_{[2]}(M)$,
i.e.~an equivariant map $(\Pi T)^{2}M\rightarrow\Gamma(T_{\vec{l}}^{*}\mathbb{R}^{0|2})$,
using a section $s$ of $T_{\vec{l}}^{*}M$: for any map $\phi:\mathbb{R}^{0|2}\rightarrow M$
we have the section $\phi^{*}s$ of $T_{\vec{l}}^{*}\mathbb{R}^{0|2}$,
i.e.~we have a map from $(\Pi T)^{2}M=M^{\mathbb{R}^{0|2}}$ to $\Gamma(T_{\vec{l}}^{*}\mathbb{R}^{0|2})$.
This is in some sense the central idea of this section.

Let us now try directly to find the space $K_{\vec{l}}(M)=Hom_{\left(\mathbb{R}^{0|2}\right)^{\mathbb{R}^{0|2}}}(\Gamma(T_{\vec{l}}^{*}\mathbb{R}^{0|2})^{*},\Omega_{[2]}(M))$.
We have already found a map $\Gamma(T_{\vec{l}}^{*}M)\rightarrow K_{\vec{l}}(M)$;
we shall construct a natural vector bundle $\widetilde{T_{\vec{l}}^{*}}M$
over $M$ such that $K_{\vec{l}}(M)=\Gamma(\widetilde{T_{\vec{l}}^{*}}M)$.
The bundle $\widetilde{T_{\vec{l}}^{*}}M$ will come with a natural
filtration $F_{1}\subset F_{2}\subset\dots\subset\widetilde{T_{\vec{l}}^{*}}M$,
such that $F_{1}=T_{\vec{l}}^{*}M$ and the quotients $F_{i}/F_{i-1}$
are $T_{\vec{k}}^{*}M$'s for various $\vec{k}$'s. The rule for getting
these $\vec{k}$'s from $\vec{l}$ should be obvious from this example:

\begin{center}\includegraphics[%
  width=6cm,
  keepaspectratio]{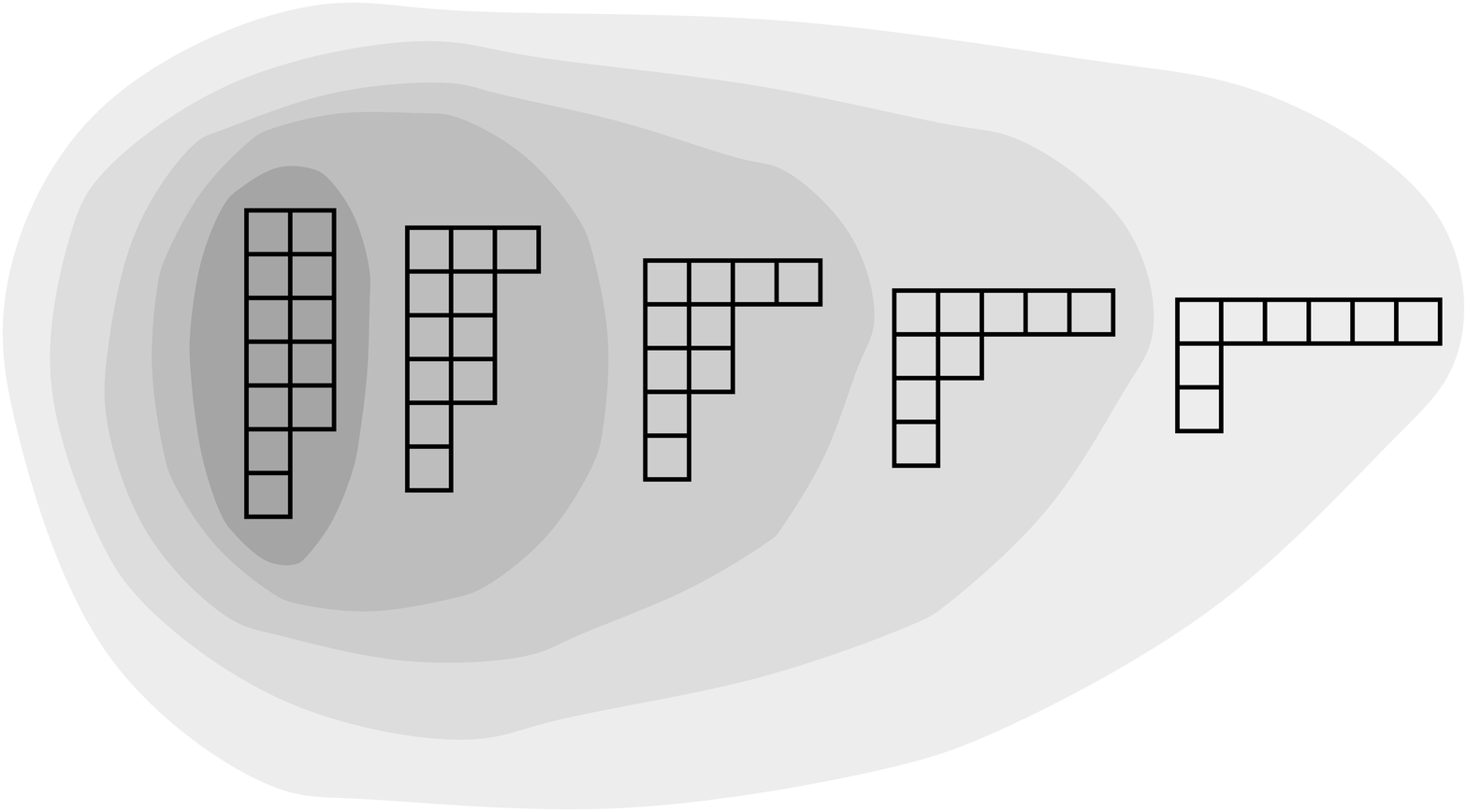}\end{center}

\noindent In other words, we start with a two-column table $\vec{l}$
and we keep removing squares from the columns and adding a square
to the first row until the second column has length one.

Let us finally find the spaces $K_{\vec{l}}(M)$, describing the bundles
$\widetilde{T_{\vec{l}}^{*}}M$ at the same time. An intertwining
map $\Gamma(T_{\vec{l}}^{*}\mathbb{R}^{0|2})^{*}\rightarrow\Omega_{[2]}(M)$
is the same as a gorm with the weight $\vec{l}^{T}$ ($\vec{l}^{T}$
denotes the Young table $\vec{l}$ reflected with respect to the diagonal,
i.e.~$\vec{l}^{T}$ is a two-row table), annulled by $E_{1}^{2}=\xi_{1}^{i}\partial_{\xi_{2}^{i}}$
(i.e.~a highest-weight gorm) and by the two operators $R_{a}=\xi_{a}^{i}\partial_{y^{i}}$.
This gorm is the image of the element of $\Gamma(T_{\vec{l}}^{*}\mathbb{R}^{0|2})^{*}$
that assigns to any section of $T_{\vec{l}}^{*}\mathbb{R}^{0|2}$
the highest-weight component of its value at the origin of $\mathbb{R}^{0|2}$.
The reason for $\vec{l}^{T}$ is that the action of $S_{N}$ on $(T^{*})^{\otimes N}\mathbb{R}^{0|2}$
is modified by the sign rule, which is equivalent to the reflection
of Young tables.

The vector bundle $\widetilde{T_{\vec{l}}^{*}}M$ is thus the vector
bundle of gorms of bidegree $\vec{l}^{T}$ that are annulled by $\xi_{1}^{i}\partial_{\xi_{2}^{i}}$
and by $R_{a}$'s. The filtration on $\widetilde{T_{\vec{l}}^{*}}M$
can be described as follows. Let us embed $M$ to $(\Pi T)^{2}M$
as the space of constant maps $\mathbb{R}^{0|2}\rightarrow M$; in
coordinates, it is given by setting $\xi$'s and $y$'s to zero. On
$\Omega_{[2]}(M)$ we have the decreasing filtration by the order
of vanishing on $M$ (i.e.~by the number of $y$'s and $\xi$'s).
The filtration on $\widetilde{T_{\vec{l}}^{*}}M$ is the restriction
of this filtration.

\subsubsection{An example: $\widetilde{T_{\boxplus}^{*}}M$}

Let us compute an example, for the Young table $\boxplus$. A general
gorm with bidegree $(2,2)$ is of the form $a_{ij}(x)\, y^{i}y^{j}+b_{ijk}(x)\, y^{i}\xi_{1}^{j}\xi_{2}^{k}+c_{ijkl}(x)\,\xi_{1}^{i}\xi_{1}^{j}\xi_{2}^{k}\xi_{2}^{l}$.
It is annulled by $E_{2}^{1}=\xi_{2}^{i}\partial_{\xi_{1}^{i}}$ iff
$b_{ijk}$ is symmetric in $jk$ and $c_{ijkl}$ has the symmetries
of the Riemann curvature tensor (i.e.~it has the symmetries given
by the Young table $\boxplus$). It is annulled by $R_{a}=\xi_{a}^{i}\partial_{y^{i}}$
iff $a_{ij}=0$ and $b_{ijk}$ is completely symmetric, i.e.~it has
the symmetries given by the Young table $\square\!\square\!\square$.
The pair $(b_{ijk},c_{ijkl})$ is a section of $\widetilde{T_{\boxplus}^{*}}M$.
The subbundle $T_{\boxplus}^{*}M\subset\widetilde{T_{\boxplus}^{*}}M$
is given by $b_{ijk}=0$; the quotient $\widetilde{T_{\boxplus}^{*}}M/T_{\boxplus}^{*}M$
is clearly $T_{\square\!\square\!\square}^{*}M$.

Similar coordinate computation can be done for arbitrary two-column
Young table $\vec{l}$; it gives the filtration on $\widetilde{T_{\vec{l}}^{*}}M$
with the Young tables as drawn on the picture above.

\subsubsection{Decomposition of the generic part}

We can conclude that the generic part of $\Omega_{[2]}(M)$ is \[
\bigoplus_{\vec{l}}\Gamma(T_{\vec{l}}^{*}\mathbb{R}^{0|2})^{*}\otimes\Gamma(\widetilde{T_{\vec{l}}^{*}}M),\]
where we sum over all two-column Young tables $\vec{l}$. To make
this formula more symmetric, let us notice that for any such $\vec{l}$,
$\widetilde{T_{\vec{l}}^{*}}\mathbb{R}^{0|2}=T_{\vec{l}}^{*}\mathbb{R}^{0|2}$,
since $T_{\vec{k}}^{*}\mathbb{R}^{0|2}=0$ for any $\vec{k}$ with
at least 3 columns. The generic part of $\Omega_{[2]}(M)$ is thus\[
\bigoplus_{\vec{l}}\Gamma(\widetilde{T_{\vec{l}}^{*}}\mathbb{R}^{0|2})^{*}\otimes\Gamma(\widetilde{T_{\vec{l}}^{*}}M).\]

\subsubsection{The non-generic part and differential forms}

Let us finish the decomposition of $\Omega_{[2]}(M)$ by describing
its non-generic part. If $\alpha\in\Omega(M)$ then for any map $\phi:\mathbb{R}^{0|2}\rightarrow M$
we have the differential form $\phi^{*}\alpha$ on $\mathbb{R}^{0|2}$,
i.e.~$\alpha$ gives us an equivariant map $(\Pi T)^{2}M\rightarrow\Omega(\mathbb{R}^{0|2})$,
i.e.~an intertwining map $\Omega(\mathbb{R}^{0|2})^{*}\rightarrow\Omega_{[2]}(M)$.
This map depends on $\alpha$, i.e.~we have found an intertwining
map \begin{equation}
\Omega(\mathbb{R}^{0|2})^{*}\otimes\Omega(M)\rightarrow\Omega_{[2]}(M).\label{eq:formy}\end{equation}
The image of this map is easily seen to be the whole non-generic part
of $\Omega_{[2]}(M)$. The map is not injective ($\Omega(\mathbb{R}^{0|2})^{*}$
is not irreducible); we have to describe its kernel.

The kernel is given by the following simple fact: the map from $\Omega(M)$
to the space of equivariant maps $(\Pi T)^{2}M\rightarrow\Omega(\mathbb{R}^{0|2})$
is a morphism of cochain complexes; in other words, the map (\ref{eq:formy})
is also a morphism of cochain complexes, where $\Omega_{[2]}(M)$
is understood as a cochain complex in the trivial way, i.e.~it is
entirely in degree 0. It means that the elements of $\Omega(\mathbb{R}^{0|2})^{*}\otimes\Omega(M)$
of non-zero degree are mapped to zero, and the same is the fate of
the exact elements of degree 0. One can easily verify that this is
the entire kernel, hence the non-generic part of $\Omega_{[2]}(M)$
is isomorphic to the degree-0 part of $\Omega(\mathbb{R}^{0|2})^{*}\otimes\Omega(M)$
modulo exact elements,\[
\frac{\bigoplus_{k}\Omega^{k}(\mathbb{R}^{0|2})^{*}\otimes\Omega^{k}(M)}{d\!\left(\bigoplus_{k}\Omega^{k}(\mathbb{R}^{0|2})^{*}\otimes\Omega^{k-1}(M)\right)}.\]

\subsubsection{The entire decomposition}

If we put the generic and the non-generic part of $\Omega_{[2]}(M)$
together, we have the isomorphism\[
\Omega_{[2]}(M)\cong\left({\textstyle \bigoplus_{\vec{l}}}\;\Gamma(T_{\vec{l}}^{*}\mathbb{R}^{0|2})^{*}\otimes\Gamma(\widetilde{T_{\vec{l}}^{*}}M)\right)\oplus\frac{\bigoplus_{k}\Omega^{k}(\mathbb{R}^{0|2})^{*}\otimes\Omega^{k}(M)}{d\!\left(\bigoplus_{k}\Omega^{k}(\mathbb{R}^{0|2})^{*}\otimes\Omega^{k-1}(M)\right)};\]
this isomorphism is $\left(\mathbb{R}^{0|2}\right)^{\mathbb{R}^{0|2}}$-equivariant
and functorial in $M$ (in particular, it is $\mathit{Diff}(M)$-equivariant).

\subsection{Derivations on gorms as a module of the crossed product of $C^{\infty}(\mathbb{R}^{0|2})$
with $\mathit{Diff}(\mathbb{R}^{0|2})$}

In this section we briefly apply the general formulas of section \ref{sub:Vector-fields-on}
to vector fields on $(\Pi T)^{2}M$. The action of $\theta_{1}$,
$\theta_{2}$, $d_{1}$ and $d_{2}$ on the space $\mathit{Der}(\Omega_{[2]}(M))$
of these vector fields generate an action of a Clifford algebra. As
a result, we have an isomorphism\[
\mathit{Der}(\Omega_{[2]}(M))\cong\mathit{Der}(\Omega_{[2]}(M))_{0}\otimes{\textstyle \bigwedge}(\mathbb{R}^{2}),\]
 where $\mathit{Der}(\Omega_{[2]}(M))_{0}$ is the space of vector
fields annulled by $\theta_{1}$ and $\theta_{2}$ and $\mathbb{R}^{2}$
is the vector space with the basis $d_{1}$ and $d_{2}$. A vector
field $A^{i}\partial_{x^{i}}+B_{a}^{i}\partial_{\xi_{a}^{i}}+C^{i}\partial_{y^{i}}$
is annulled by both $\theta$'s iff $A^{i}=B_{a}^{i}=0$. We can thus
naturally identify $\mathit{Der}(\Omega_{[2]}(M))_{0}$ with the space
of gorm-valued vector fields on $M$. Similar result hods for worms
of arbitrary level.

To give a decomposition of $\mathit{Der}(\Omega_{[2]}(M))$ as a module
of the entire crossed product of $C^{\infty}(\mathbb{R}^{0|2})$ with
$\mathit{Diff}(\mathbb{R}^{0|2})$ we would have to take into account
the action of $\mathit{Mat}(2)$ (which is easy) and of $R_{a}$'s;
the result seems to be more complicated than interesting, so we shall
not write it here.

\section{Integration and Euler characteristic}

This section is devoted to the proof of the following theorem: if
$\gamma$ is a pseudodifferential gorm on a connected manifold $M$
(i.e.~a smooth function on $(\Pi T)^{2}M$) such that $d_{1}\gamma=d_{2}\gamma=0$
(this clearly implies that the restriction of $\gamma$ to $M\subset(\Pi T)^{2}M$,
$\gamma|_{M}$, is a constant), and if moreover $\gamma$ is integrable
and $M$ compact, then\begin{equation}
\int\gamma=\frac{\gamma|_{M}}{2}(-\pi)^{m/2}S_{m}\,\chi(M),\label{eq:integral}\end{equation}
where $\chi(M)$ is the Euler characteristic of $M$, $m$ is the
dimension of $M$ and $S_{m}$ is the area of the unit $m$-dimensional
sphere.

Let us start with a special case. Let $b_{ij}(x)$ be a Riemann metric
on $M$ and let $\beta=b_{ij}(x)\, d_{1}x^{i}\, d_{2}x^{j}$; we will
prove the theorem for $\gamma=e^{d_{1}d_{2}\beta}$. If we choose
local coordinates so that $b_{ij,k}=0$ at a given point (Riemann
normal coordinates would do), then a simple computation gives that
at that point \[
d_{1}d_{2}\beta=-b_{ij}\, d_{1}d_{2}x^{i}\, d_{1}d_{2}x^{j}-\frac{1}{2}R_{ijkl}\, d_{1}x^{i}\, d_{1}x^{j}\, d_{2}x^{k}\, d_{2}x^{l}\]
where $R_{ijkl}$ is the curvature of $b_{ij}$.

If $M$ is compact then $e^{d_{1}d_{2}\beta}$ is integrable. To compute
the integral, first pass to Riemann normal coordinates at a point
and integrate over $y$'s and $\xi$'s; we end up with the Pfaffian
of the curvature, whose integral is well known to be a multiple of
the Euler characteristic $\chi(M)$ of $M$. The result is really\begin{equation}
\int e^{d_{1}d_{2}\beta}=\frac{1}{2}(-\pi)^{m/2}S_{m}\,\chi(M).\label{eq:metricEuler}\end{equation}

It is easy to prove directly that the integral (\ref{eq:metricEuler})
is independent of the choice of the metric. If we add to $\beta$
an infinitesimal $\alpha$ then $e^{d_{1}d_{2}(\beta+\alpha)}-e^{d_{1}d_{2}\beta}=e^{d_{1}d_{2}\beta}d_{1}d_{2}\alpha=d_{1}(e^{d_{1}d_{2}\beta}d_{2}\alpha)$
and $\int d_{1}(e^{d_{1}d_{2}\beta}d_{2}\alpha)=0$ since $\int u^{\flat}\gamma=0$
for any $u$ and any integrable $\gamma$ (see section \ref{sub:Integration-of-gorms}).

To prove (\ref{eq:integral}) generally, we first have to prove that
if $\delta$ grows (say) at most polynomially in $y$'s, $d_{1}\delta=d_{2}\delta=0$
and $\delta|_{M}=0$, then $\int e^{d_{1}d_{2}\beta}\delta=0$. Indeed,
since $\delta|_{M}=0$, we can find an $\epsilon$ such that $\delta=E\epsilon$,
where $E=E_{1}^{1}+E_{2}^{2}$ is the generator of the scaling $(x^{i},\,\xi_{a}^{i},\, y^{i})\mapsto(x^{i},\,\lambda\xi_{a}^{i},\,\lambda^{2}y^{i})$.
Since $E=-[d_{1},R_{1}]-[d_{2},R_{2}]$, $\delta=-d_{1}R_{1}\epsilon-d_{2}R_{2}\epsilon$
and $\int e^{d_{1}d_{2}\beta}\delta=\int-d_{1}(e^{d_{1}d_{2}\beta}R_{1}\epsilon)-d_{2}(e^{d_{1}d_{2}\beta}R_{2}\epsilon)=0$.

Finally, to prove (\ref{eq:integral}) we just set $\delta=\gamma-\gamma|_{M}$,
multiply $\beta$ by a constant $s>0$ and take the limit $s\rightarrow0_{+}$.

\section{Beyond homological algebra}

Differential forms certainly play an important role in topology; they
are the basic example of homological algebra, and also of its non-linear
generalizations (consider e.g.~Maurer-Cartan equation for flat connections,
or Sullivan's rational homotopy theory). Since in this paper we {}``explained''
and generalized differential forms, it is also natural to {}``explain''
and generalize homological algebra and its non-linear analogs. As
we will see there is a closely related problem: to define differential
forms (and gorms and worms) for some generalized manifolds, namely
for contravariant functors and more generally for stacks.

It shouldn't be surprising that homological algebra is closely
connected with supermanifolds with a right action of the
supersemigroup $\left(\mathbb{R}^{0|1}\right)^{\mathbb{R}^{0|1}}$
(these objects are really taken from \cite{Sul}); we shall discuss
this connection in section \ref{sub:Example:--n=3D1}. The
generalization is simply to replace $\mathbb{R}^{0|1}$ with
$\mathbb{R}^{0|n}$. The description of differential forms as
functions on $M^{\mathbb{R}^{0|1}}$ gives a new point of view on
(a part of) homological algebra via representability of functors.
It is closely connected with the problem of differential forms
(and gorms and worms) on contravariant functors; we shall discuss
it in section \ref{sub:Representability-of-functors}.

Let us introduce a part of the picture. By $\mathcal{S}_{[n]}$ we
denote the category of supermanifolds with right action of the
supersemigroup $\left(\mathbb{R}^{0|n}\right)^{\mathbb{R}^{0|n}}$;
morphisms in $\mathcal{S}_{[n]}$ are equivariant maps. We shall
view the objects of $\mathcal{S}_{[n]}$ as generalized
supermanifolds; this idea is taken directly from Sullivan's
rational homotopy theory \cite{Sul} (generalized {\em manifolds}
are then such generalized supermanifolds, on which the parity
involution $\mathbb{R}^{0|n}\rightarrow\mathbb{R}^{0|n}$ acts as
the parity involution). The idea is as follows. The functor
$X\mapsto X^{\mathbb{R}^{0|n}}$ is a fully faithful embedding of
the category $\mathcal{S}$ of supermanifolds to the category
$\mathcal{S}_{[n]}$. The category $\mathcal{S}_{[n]}$ thus becomes
an extension of $\mathcal{S}$ and we can regard objects of
$\mathcal{S}_{[n]}$ as generalized supermanifolds; true
supermanifolds are the objects of $\mathcal{S}_{[n]}$ of the form
$(\Pi T)^{n}X=X^{\mathbb{R}^{0|n}}$. In other words, we shall
treat the objects of $\mathcal{S}_{[n]}$ as if they were of the
form $X^{\mathbb{R}^{0|n}}$ for some $X$. For example, if we have
two objects $\mathsf{X},\mathsf{Y}$ of $\mathcal{S}_{[n]}$ and two
equivariant maps (i.e.~two morphisms) between them, a
\emph{homotopy} between the maps is an equivariant map
$\mathsf{X}\times(\Pi T)^{n}I\rightarrow\mathsf{Y}$ that restricts
to the two maps at the endpoints of $I$.

This idea will become more complete and convincing in section \ref{sub:Representability-of-functors}
when we identify objects of $\mathcal{S}_{[n]}$ with the functors
$\mathcal{S}^{op}\rightarrow\mathcal{S}$ they represent. In particular,
we will get a whole chain of fully faithful embeddings\[
\mathcal{S}\rightarrow\mathcal{S}_{[1]}\rightarrow\mathcal{S}_{[2]}\rightarrow\mathcal{S}_{[3]}\rightarrow\cdots.\]
In section \ref{sub:stacks} we will extend it to Lie groupoids and
to the stacks they represent.

\subsection{Example: $n=1$ (the case of homological algebra)\label{sub:Example:--n=3D1}}

Let us consider the case of $n=1$, corresponding to differential
forms, homological algebra etc. Let us start with linear actions of
$\left(\mathbb{R}^{0|1}\right)^{\mathbb{R}^{0|1}}\!$. A right representation
of $\left(\mathbb{R}^{0|1}\right)^{\mathbb{R}^{0|1}}\!\!\!$, i.e.~a
{}``generalized vector space'', is the same as a non-negatively
graded chain complex ($\theta\partial_{\theta}$ is the degree and
$\partial_{\theta}$ is the differential). Two linear equivariant
maps $\mathsf{V}\rightrightarrows\mathsf{W}$ are linearly homotopic,
i.e.~they are connected by an equivariant map $\mathsf{V}\times\Pi TI\rightarrow\mathsf{W}$
that is a linear map $\mathsf{V}\rightarrow\mathsf{W}$ parametrized
by $\Pi TI$, iff the two morphisms of chain complexes are homotopic
in the usual algebraic sense. If we define homotopy groups of $\mathsf{V}$
using equivariant maps $\Pi TS^{k}\rightarrow\mathsf{V}$, they turn
out to be the homology groups of the complex $\mathsf{V}$.

Let us now pass to some non-linear actions of $\left(\mathbb{R}^{0|1}\right)^{\mathbb{R}^{0|1}}\!$.
If $\mathfrak{g}$ is a Lie algebra then on $\Pi\mathfrak{g}$ we
have a canonical action of $\left(\mathbb{R}^{0|1}\right)^{\mathbb{R}^{0|1}}\!\!\!$,
given by identification of $\Pi\mathfrak{g}$ with $(\Pi TG)/G$ (in
other words, $C^{\infty}(\Pi\mathfrak{g})=\bigwedge\mathfrak{g}^{*}$,
$\theta\partial_{\theta}$ acts as the degree and $\partial_{\theta}$
as the Chevalley-Eilenberg differential). An equivariant map $\Pi TM\rightarrow\Pi\mathfrak{g}$
is the same as a flat $\mathfrak{g}$-connection on $M$: any map
$\Pi TM\rightarrow\Pi\mathfrak{g}$ is a $\mathfrak{g}$-valued differential
form on $M$, and equivariance is easily seen to express the fact
that it is a 1-form satisfying the Maurer-Cartan equation. The fundamental
group of $\Pi\mathfrak{g}$ is therefore the 1-connected Lie group
$G$, and its higher homotopy groups are the higher homotopy groups
of $G$.

As a little generalization, if $A\rightarrow N$ is a Lie algebroid
then again $\Pi A$ is an object of $\mathcal{S}_{[1]}$ (this was
observed by Vaintrob \cite{vain}). An equivariant map $\Pi TM\rightarrow\Pi A$
is the same as a Lie algebroid morphism $TM\rightarrow A$; the fundamental
groupoid of $\Pi A$ is therefore the corresponding Lie groupoid $\Gamma$
with 1-connected fibres (if it exists).

For more details and examples with interesting higher homotopies,
see \cite{se2}.

\subsection{Representability of functors and their approximations\label{sub:Representability-of-functors}}

For any supermanifolds $X$, $Y$ we have the supermanifold $Y^{X}$
of all maps $X\rightarrow Y$ (we shall ignore all problems connected
with the fact that $Y^{X}$ is almost always infinite-dimensional,
as they are inessential for our purposes; just imagine that we designed
our category $\mathcal{S}$ of supermanifolds so that it contains
$Y^{X}$ for any objects $X$ and $Y$). Hence any object $Y\in\mathcal{S}$
gives us a functor $\hat{Y}:\mathcal{S}^{op}\rightarrow\mathcal{S}$
given by \[
\hat{Y}(X)=Y^{X}.\]
Functors of this form are called representable. $Y$ can be reconstructed
as $\hat{Y}(\mathit{point})$, and for example $(\Pi T)^{n}Y$ as
$\hat{Y}(\mathbb{R}^{0|n})$.

Let $\mathcal{SS}$ denote the category of all functors $F:\mathcal{S}^{op}\rightarrow\mathcal{S}$.%
\footnote{We understand these contravariant functors in the strong sense: for
any two $X,Y\in\mathcal{S}$ we have a map of supermanifolds $Y^{X}\times F(Y)\rightarrow F(X)$.
Equivalently (using parametrized maps instead of map spaces), for
any triple $X,Y,Z$ and any map $Z\times X\rightarrow Y$ we have
a map $Z\times F(Y)\rightarrow F(X)$. Morphisms between functors
are understood in the strong sense too: a morphism $F_{1}\rightarrow F_{2}$
is a morphism $F_{1}(X)\rightarrow F_{2}(X)$ for each $X$ such that
for any $Z\times X\rightarrow Y$ the square\[
\begin{array}{ccc}
F_{1}(X) & \rightarrow & F_{2}(X)\\
\uparrow &  & \uparrow\\
Z\times F_{1}(Y) & \rightarrow & Z\times F_{2}(Y)\end{array}\]
commutes.%
} We have described a functor $\mathcal{S}\rightarrow\mathcal{SS}$,
$Y\mapsto\hat{Y}$; by Yoneda lemma it is a fully faithful embedding.
We can (and will) identify objects of $\mathcal{S}$ with the functors
they represent. It is a standard idea to view objects of $\mathcal{SS}$
as generalized objects of $\mathcal{S}$ (a generalized \emph{manifold}
is a functor $F:\mathcal{S}^{op}\rightarrow\mathcal{S}$ that preserves
the parity involution). In other words, we will understand $F(X)$
as the space of maps from $X$ to some generalized space corresponding
to $F$.%
\footnote{$F(X)$ can be seen as just an approximation to the generalized space
of maps $F^{X}$ defined by $F^{X}(Y)=F(X\times Y)$, i.e.~$F(X)=F^{X}(\mathit{point})$%
} From this point of view we should define level-$n$ worms on $F$
as functions on $F(\mathbb{R}^{0|n})$. Notice that $\left(\mathbb{R}^{0|n}\right)^{\mathbb{R}^{0|n}}$
acts on $F(\mathbb{R}^{0|n})$ from the right, i.e.~$F(\mathbb{R}^{0|n})$
is an object of $\mathcal{S}_{[n]}$. In fact the semigroup $\left(\mathbb{R}^{0|n}\right)^{\mathbb{R}^{0|n}}$
can be understood as the full subcategory of $\mathcal{S}$ with just
one object, $\mathbb{R}^{0|n}$, and $\mathcal{S}_{[n]}$ as the category
of contravariant functors from this subcategory to $\mathcal{S}$;
we restricted $F$ to this subcategory.

Although we understand $F(\mathit{point})$ as the space of points
of $F$, the functor $F$ is not uniquely specified by $F(\mathit{point})$
(otherwise all functors would have to be representable), nor is it
specified by $F(\mathbb{R}^{0|n})$ for any $n$. Nevertheless we
can use $F(\mathbb{R}^{0|n})$'s to approximate $F$, since for any
$X$ we have the map \[
X^{\mathbb{R}^{0|n}}\times F(X)\rightarrow F(\mathbb{R}^{0|n}),\]
i.e.~a map from $F(X)$ to the superspace of $\left(\mathbb{R}^{0|n}\right)^{\mathbb{R}^{0|n}}$-equivariant
maps $X^{\mathbb{R}^{0|n}}\rightarrow F(\mathbb{R}^{0|n})$,\begin{equation}
F(X)\rightarrow\mathit{Hom}_{\mathcal{S}_{[n]}}(X^{\mathbb{R}^{0|n}},F(\mathbb{R}^{0|n})).\label{eq:approx}\end{equation}

Now we can state the definitions. For any object
$\mathsf{Y}\in\mathcal{S}_{[n]}$, let
$\hat{\mathsf{Y}}\in\mathcal{SS}$ (the \emph{functor represented
by $\mathsf{Y}$}) be the functor given by $\hat{\mathsf{Y}}(X)=$
the superspace of all
$\left(\mathbb{R}^{0|n}\right)^{\mathbb{R}^{0|n}}$-equivariant
maps from $X^{\mathbb{R}^{0|n}}$ to $\mathsf{Y}$,\begin{equation}
\hat{\mathsf{Y}}(X)=\mathit{Hom}_{\mathcal{S}_{[n]}}(X^{\mathbb{R}^{0|n}},\mathsf{Y}).\label{eq:representable}\end{equation}
(the functor $\mathcal{S}_{[n]}\rightarrow\mathcal{SS}$,
$\mathsf{Y}\mapsto\hat{\mathsf{Y}}$, is the right adjoint of the
restriction functor $\mathcal{SS}\rightarrow\mathcal{S}_{[n]}$).
Functors of the form $\hat{\mathsf{Y}}$ will be called
\emph{representable at level n}. Notice that they can be expressed
in terms of worms of level $n$ (if we choose coordinates on
$\mathsf{Y}$, a map $X^{\mathbb{R}^{0|n}}\rightarrow\mathsf{Y}$
becomes a collection of functions on $X^{\mathbb{R}^{0|n}}$,
i.e.~a collection of level-$n$ worms on $X$). Representability at
level 0 is, of course, the ordinary representability. If
$F:\mathcal{S}^{op}\rightarrow\mathcal{S}$ is representable at
level $n$, the corresponding object
$\mathsf{Y}\in\mathcal{S}_{[n]}$ can be found as
$F(\mathbb{R}^{0|n})$. For any functor
$F:\mathcal{S}^{op}\rightarrow\mathcal{S}$, the functor
$F_{[n]}=\widehat{F(\mathbb{R}^{0|n})}$ will be called \emph{the
n-th approximation of F.} The equation (\ref{eq:approx}) gives us
a natural morphism $F\rightarrow F_{[n]}$; $F$ is representable at
level $n$ iff the morphism is an isomorphism.

The morphism $F\rightarrow F_{[n]}$ becomes an isomorphism when we
restrict $F$ and $F_{[n]}$ to the full subcategory
$\mathcal{D}_{n}\subset\mathcal{S}$ of supermanifolds of dimension
at most $0|n$ (we have to show that
$F(\mathbb{R}^{0|k})\rightarrow F_{[n]}(\mathbb{R}^{0|k})$ is a
diffeomorphism whenever $k\leq n$; for $k=n$ it is tautological,
and for other $k$'s it is enough to choose maps
$\mathbb{R}^{0|k}\rightarrow\mathbb{R}^{0|n}\rightarrow\mathbb{R}^{0|k}$
that compose to identity on $\mathbb{R}^{0|k}$). As a consequence,
representability at level $n$ implies representability at all
higher levels. By taking successive approximations we get a chain
of morphisms \begin{equation} \cdots\rightarrow F_{[3]}\rightarrow
F_{[2]}\rightarrow F_{[1]}\rightarrow
F_{[0]};\label{eq:mchain}\end{equation} together with the
morphisms $F\rightarrow F_{[n]}$ it forms a commutative diagram.

As a final remark, we could repeat these definitions from a more
natural point of view. The idea is to approximate functors
$\mathcal{S}^{op}\rightarrow\mathcal{S}$ (objects of
$\mathcal{SS}$) by their restrictions to the full subcategories
$\mathcal{D}_{n}$ of $\mathcal{S}$ defined above that form a chain
of inclusions\[
\mathcal{D}_{0}\subset\mathcal{D}_{1}\subset\mathcal{D}_{2}\subset\mathcal{D}_{3}\subset\cdots\subset\mathcal{S}.\]
The category $\mathcal{D}_{n}\mathcal{S}$ of all functors
$\mathcal{D}_{n}^{op}\rightarrow\mathcal{S}$ is equivalent to
$\mathcal{S}_{[n]}$ (the equivalence
$\mathcal{D}_{n}\mathcal{S}\rightarrow\mathcal{S}_{[n]}$ is given
by restriction), so we would get equivalent definitions.

\subsection{Examples of approximations\label{sub:Examples-of-approximations}}

Let $G$ be a Lie group and let $F:\mathcal{S}^{op}\rightarrow\mathcal{S}$
be given by $F(X)=(G^{X})/G.$ To compute the $n$-th approximation
of $F$ we just have to compute the space $F(\mathbb{R}^{0|n})$ and
the right action of $\left(\mathbb{R}^{0|n}\right)^{\mathbb{R}^{0|n}}$
on this space. Since $F(\mathit{point})=\mathit{point}$, $F_{[0]}(X)=\mathit{point}$.
The first approximation is more interesting: $F(\mathbb{R}^{0|1})=(\Pi TG)/G=\Pi\frak{g}$,
therefore $F_{[1]}(X)$ is the space of flat $\frak{g}$-connections
on $X$. Since locally (in $X$) one cannot distinguish $F$ from
$F_{[1]}$, all higher approximations of $F$ are just $F_{[1]}$.

As a small generalization, let $\Gamma$ be a Lie groupoid. For any
$X$ let $X\times X$ be the pair groupoid (with $X$ as the space
of objects and with one arrow between any two objects) and finally
let $F(X)=\mathit{Hom}(\mathit{X\times X},\Gamma).$ Then $F(\mathit{point})$
is the base of $\Gamma$ (the space of its objects), hence $F_{[0]}$
is just the functor represented by the base. To compute $F_{[1]}$
notice that $F(\mathbb{R}^{0|1})=\Pi A$ where $A$ is the Lie algebroid
corresponding to $\Gamma$. Thus we found that $F_{[1]}(X)$ is the
space of all Lie algebroid morphisms $TX\rightarrow A$. Since locally
we cannot distinguish between these Lie algebroid morphisms and Lie
groupoid morphisms $X\times X\rightarrow\Gamma$, all higher approximations
of $F$ are again equal to $F_{[1]}$.

The next example is trivially representable at level 1, but it is
interesting for other reasons. Let $G$ be a Lie group,
$\mathfrak{g}$ its Lie algebra, and let $F(X)$ be the space of
$\mathfrak{g}$-connections on $X$ (i.e.~the space of
$\mathfrak{g}$-valued 1-forms on $X$). On $F(X)$ we have action of
the group $G^{X}$ (by gauge transformations); if we understand $F$
as a generalized space, it means that $G$ acts on $F$. The algebra
of differential forms on $F$, i.e.~of functions on
$F(\mathbb{R}^{0|1})$, is the Weil algebra $W(\mathfrak{g})$; the
actions of $\left(\mathbb{R}^{0|1}\right)^{\mathbb{R}^{0|1}}\!$
and of $\Pi TG=G^{\mathbb{R}^{0|1}}$ give rise to its standard
$G$-differential algebra structure. If $M$ is a manifold with an
action of $G$, we can consider the generalized space $F\times M$
(given by $(F\times M)(X)=F(X)\times M^{X}$); since $G$ acts on
both $F$ and $M$, it acts also on $F\times M$. The complex of
basic forms on $F\times M$, i.e.~of $\Pi TG$-invariant functions
on $(F\times M)(\mathbb{R}^{0|1})=F(\mathbb{R}^{0|1})\times\Pi
TM$, is the basic subcomplex in the Weil model of equivariant
cohomology (notice that $F$ behaves as if it were $EG$). To get
Cartan model, notice that any connection on $\mathbb{R}^{0|1}$ can
be made to vanish at the origin, using a suitable gauge
transformation; the space of such connections can be identified
with $\mathfrak{g}$ (any such connection is of the form $t\theta
d\theta$, where $t\in\mathfrak{g}$). After we impose this
condition, only the group of constant gauge transformations
$G\subset\Pi TG$ remains. $\Pi TG$-invariant functions on
$F(\mathbb{R}^{0|1})\times\Pi TM$ can thus be identified with
$G$-invariant functions on $\mathfrak{g}\times\Pi TM$; the latter
is the Cartan model.

Let us pass to some examples where $F_{[2]}$ is different from $F_{[1]}$.
Let $F(X)=\Gamma(S^{2}T^{*}X)$ (or $\Gamma(T_{\square\!\square}^{*}X)$
in the notation of section \ref{sec:representation}). Then $F(\mathit{point})=F(\mathbb{R}^{0|1})=0$.
On the other hand, as we have found in section \ref{sec:representation},
$F=F_{[2]}$ since $F_{[2]}(X)=\Gamma(\widetilde{T_{\square\!\square}^{*}}X)$
and $\widetilde{T_{\square\!\square}^{*}}=T_{\square\!\square}^{*}$.
The functor $F$ is thus representable at level 2. As another example,
let $F(X)=\Gamma(T_{\boxplus}^{*}X)$. Then again $F(\mathit{point})=F(\mathbb{R}^{0|1})=0$,
but $F_{[2]}(X)=\Gamma(\widetilde{T_{\boxplus}^{*}}X)$. This time
$F$ is different from its second approximation; one can prove that
it is representable at level 3.

Let us finish with a simple example of a functor $F$ for which the
chain of morphisms (\ref{eq:mchain}) doesn't stablilize. Let
$F(X)=C^\infty(X\times X)$. Then $F_{[k]}(X)=\Gamma(J^k(X))$,
where $J^k(X)\rightarrow X$ is the vector bundle of $k$-jets of
functions on $X$ (i.e.~$F_{[k]}(X)$ is the space of functions on
the $k$th formal neighbourhood of the diagonal in $X\times X$.)

\subsection{Approximations and representability of stacks \textcolor{red}{\label{sub:stacks}}}

This section can be seen in two ways. Above we defined the categories
$\mathcal{S}_{[n]}$ as generalizations of the category $\mathcal{S}$
of supermanifolds. Here we extend these generalizations from supermanifolds
to Lie supergroupoids (any (super)manifold can be seen as a Lie (super)groupoid,
with identity arrows only). We define 2-categories $\mathcal{G}_{[n]}$
as generalizations of the 2-category $\mathcal{G}$ of Lie supergroupoids.
The objects of $\mathcal{G}_{[n]}$ will be Lie supergroupoids over
the supersemigroup $\left(\mathbb{R}^{0|n}\right)^{\mathbb{R}^{0|n}}$
(roughly speaking, Lie supergroupoids on which $\left(\mathbb{R}^{0|n}\right)^{\mathbb{R}^{0|n}}$
acts up to natural transformations). In the case of $n=1$ they will
include some interesting know examples where {}``$d^{2}=0$ up to
gauge transformations'', e.g.~Cartan model of equivariant cohomology,
or quasi-Poisson groupoids.

The other point of view (extending the section \ref{sub:Representability-of-functors})
is to look at the stacks represented by the objects of $\mathcal{G}_{[n]}$.
Principal $G$-bundles (for some fixed Lie group $G$) are perhaps
the simplest example of a stack, and this stack is representable in
the appropriate way: a principal bundle $P\rightarrow M$ is the same
as a 1-morphism of groupoids $M\rightarrow G$ in the sense of Hilsum
and Skandalis (see below). We can get other interesting stacks by
considering 1-morphisms $M\rightarrow\Gamma$ for a Lie groupoid $\Gamma$.

The stack of principal $G$-bundles with a choice of a connection
is, however, not representable in this sense. Fortunately, it is representable
at level 1, i.e.~by an object of $\mathcal{G}_{[1]}$. This object
can be found in a tautological way (just as in \ref{sub:Representability-of-functors},
where we would find an object of $\mathcal{S}_{[1]}$, corresponding
to a functor $F$, as $F(\mathbb{R}^{0|1})$), as the supergroupoid
of all principal $G$-bundles over $\mathbb{R}^{0|1}$ with a choice
of connection. Any principal $G$-bundle over $\mathbb{R}^{0|1}$
is trivializable so we can consider just connections on the trivial
bundle (and get equivalent groupoid). The objects of this groupoid
are thus $\mathfrak{g}$-valued 1-forms on $\mathbb{R}^{0|1}$, the
arrows are given by gauge transformations (i.e.~by the action of
the supergroup $\Pi TG=G^{\mathbb{R}^{0|1}}$); the supersemigroup
$\left(\mathbb{R}^{0|1}\right)^{\mathbb{R}^{0|1}}$ acts from the
right on this supergroupoid, so we get an object of $\mathcal{G}_{[1]}$.
We can get a yet smaller equivalent groupoid by considering only those
$\mathfrak{g}$-valued 1-forms that vanish at the origin (see the
example in the section \ref{sub:Examples-of-approximations}); this
is still an object of $\mathcal{G}_{[1]}$, as $\left(\mathbb{R}^{0|1}\right)^{\mathbb{R}^{0|1}}$
acts on it up to gauge transformations.

Finally we have to say that this section is not logically complete
(roughly because we do not give a definition of a {}``superstack'',
so it is not clear, how to restrict a stack to $\mathbb{R}^{0|n}$
to get an object of $\mathcal{G}_{[n]}$ (even though it is clear
in examples); on the other hand, the definition of a stack represented
by an object of $\mathcal{G}_{[n]}$ is all right). We hope that this
fault can be excused by interesting examples.

\subsubsection{Groupoids over a category and Hilsum-Skandalis morphisms}

Let us recall from \cite{SGA} that a functor $F:\mathcal{E}\rightarrow\mathcal{F}$
is a \emph{fibration of $\mathcal{F}$ by groupoids}, (or shortly,
a \emph{groupoid over $\mathcal{F}$}), if it satisfies the following
lifting property for morphisms: for any morphism $f:X\rightarrow Y$
in $\mathcal{F}$ and any object $Q$ in $\mathcal{E}$ such that
$F(Q)=Y$ there is a morphism $\tilde{f}:P\rightarrow Q$ in $\mathcal{E}$
such that $F(\tilde{f})=f$, and moreover $\tilde{f}$ is essentially
unique, i.e.~if $\tilde{f}':P'\rightarrow Q$ is another morphism
such that $F(\tilde{f}')=f$ then there is unique $h:P'\rightarrow P$
such that $f\circ h=\tilde{f}$ and $F(h)=\mathit{id}_{X}$. As a
simple example, we can take $\mathcal{E}$ to be the category of principal
$G$-bundles (with equivariant maps as morphisms) and $\mathcal{F}$
the category of manifolds; the functor $F:\mathcal{E}\rightarrow\mathcal{F}$
assigns to a principal bundle its base.

For any object $X$ of $\mathcal{F}$ let $F^{-1}(X)$ be the fibre
above $X$, i.e.~the subcategory of $\mathcal{E}$ of objects $P$
such that $F(P)=X$ and morphisms $f$ such that $F(f)=\mathit{id}_{X}$.
It is easy to see that all $F^{-1}(X)$'s are groupoids. If $f:X\rightarrow Y$
is a morphism then by lifting $f$ at all objects of $F^{-1}(Y)$
we get a functor $F^{-1}(Y)\rightarrow F^{-1}(X)$; if we choose the
lifts differently, we get an isomorphic functor. By lifting all the
morphisms of $\mathcal{F}$ we get a lax functor from $\mathcal{F}^{op}$
to the category of groupoids; in the opposite direction, any (lax)
functor from $\mathcal{F}^{op}$ to the category of groupoids gives
us a groupoid over $\mathcal{F}$.

If $\Delta_{1}$ is the {}``segment'' category with objects 0 and
1 and with only one non-identity morphism, $0\rightarrow1$, then
a groupoid over $\Delta_{1}$ is called a \emph{Hilsum-Skandalis (HS)
morphism} from $F^{-1}(1)$ to $F^{-1}(0)$; it is essentially a functor
$F^{-1}(1)\rightarrow F^{-1}(0)$, but the functor depends (up to
natural transformations) on the choice of lifts.

Groupoids form a (weak) 2-category, with HS morphisms as 1-morphisms
(if $\mathcal{C}\rightarrow\Delta_{1}$ and $\mathcal{D}\rightarrow\Delta_{1}$
are two HS morphisms from $\Gamma_{1}$ to $\Gamma_{2}$, a 2-morphism
between them is an isomorphism $\mathcal{C}\cong\mathcal{D}$ that
is identity on both $\Gamma_{1}$ and $\Gamma_{2}$). If $\Delta_{2}$
is the {}``triangle'' category, with three objects 0, 1 and 2 ,
and with three morphisms (except for identities)\[
\begin{array}{ccc}
0 & \rightarrow & 1\\
 & \searrow & \downarrow\\
 &  & 2\end{array},\]
a groupoid over $\Delta_{2}$ gives us 3 HS morphisms, $F^{-1}(1)\rightarrow F^{-1}(0)$,
$F^{-1}(2)\rightarrow F^{-1}(1)$ and $F^{-1}(2)\rightarrow F^{-1}(0)$.
The HS morphism $F^{-1}(2)\rightarrow F^{-1}(0)$ is then \emph{a
composition} of $F^{-1}(2)\rightarrow F^{-1}(1)$ and $F^{-1}(1)\rightarrow F^{-1}(0)$.
It is easy to see that a composition of two HS morphisms $\Gamma_{2}\rightarrow\Gamma_{1}\rightarrow\Gamma_{0}$
always exists and it is unique up to a canonical isomorphism: we have
to define the arrows over $0\rightarrow2$ (i.e.~arrows $X\rightarrow Z$,
where $X$ is an object of $\Gamma_{0}$ and $Z$ of $\Gamma_{2}$);
these will be, by definition, pairs of arrows $X\rightarrow Y\rightarrow Z$,
where we identify $X\rightarrow Y_{1}\rightarrow Z$ with $X\rightarrow Y_{2}\rightarrow Z$
whenever there is an arrow $Y_{1}\rightarrow Y_{2}$ such that the
two triangles\[
\begin{array}{ccccc}
X & \rightarrow & Y_{1}\\
 & \searrow & \downarrow & \searrow\\
 &  & Y_{2} & \rightarrow & Z\end{array}\]
are commutative.

\subsubsection{The 2-category of Lie groupoids and the stacks they represent }

A \emph{Lie groupoid over a category} $\mathcal{F}$ is a groupoid
over $\mathcal{F}$, $\mathcal{E}\rightarrow\mathcal{F}$, such that
for any arrow $f:X\rightarrow Y$ of $\mathcal{F}$ the set of arrows
of \emph{$\mathcal{E}$} over $f$ is a manifold, and the composition
of arrows in $\mathcal{E}$ is smooth (i.e.~the fibres of $\mathcal{E}\rightarrow\mathcal{F}$
are Lie groupoids and they act smoothly on the manifolds of arrows,
for which the composition is defined). A Hilsum-Skandalis morphism
of Lie groupoids is then simply a Lie groupoid over $\Delta_{1}$.
For example, a HS morphism $M\rightarrow G$, where $M$ is a manifold
(understood as a Lie groupoid with identity arrows only) and $G$
a Lie group, is the same as a principal $G$-bundle over $M$ (the
bundle is the manifold of arrows over $0\rightarrow1$). Lie supergroupoids
with HS morphisms form a (weak) 2-category, denoted $\mathcal{G}$.

As a generalization of principal $G$-bundles, any Lie groupoid $\Gamma$
defines a stack (a groupoid over the category of manifolds, satisfying
a sheaf-like condition): The objects of this category are HS morphisms
$M\rightarrow\Gamma$ (where $M$ runs over all manifolds) and morphisms
are commutative triangles\[
\begin{array}{ccc}
M_{1} & \rightarrow & \Gamma\\
\downarrow & \nearrow\\
M_{2}\end{array}.\]
Stacks of this form will be called \emph{representable at level} 0.
If $M$ is contractible, the groupoid of all HS morphisms $M\rightarrow\Gamma$
(i.e.~the fibre of the stack over $M$) is equivalent to the groupoid
of all (ordinary) maps $M\rightarrow\Gamma$. Notice that $\Gamma$
is equivalent to the fibre over $M=\mathit{point}$. For a review
of the 2-category $\mathcal{G}$ and of the stacks represented by
Lie groupoids, see \cite{met}.

Notice that to describe a HS morphism $M\rightarrow\Gamma$ we need
to give the space $P$ of arrows over $0\rightarrow1$, the submersion
$P\rightarrow M$ (the map sending arrows to their heads) and the
action of $\Gamma$ on $P$ (composition of arrows); we get a HS morphism
iff the submersion $P\rightarrow M$ is surjective and $\Gamma$ acts
freely and transitively on each of its fibres. A space $P$ with these
properties is (for obvious reasons) often called \emph{a principal
$\Gamma$ bundle}.

As an example, if a group $G$ acts on a space $V$, we can consider
the stack whose objects are principal $G$-bundles with equivariant
maps to $V$; this stack can be represented by the action groupoid
of $G$ on $V$.

\subsubsection{Generalized Lie groupoids and the stacks they represent}

We shall need Lie groupoids over some smooth categories. Rather than
defining general Lie categories (where, for our purposes, dimension
should be allowed to be different on different components), we do
it just in special cases that we really need, where the objects of
the base category $\mathcal{F}$ form a discrete set. So suppose $\mathcal{F}$
is a Lie category of this kind, i.e.~simply a category enriched in
the category of manifolds.%
\footnote{this means that $\mathit{Hom}(X,Y)$ is a manifold for
any objects $X$ and $Y$ of $\mathcal{F}$ and that
$\mathit{Hom}(X,Y)\times\mathit{Hom}(Y,Z)\rightarrow\mathit{Hom}(X,Z)$
is a smooth map%
} A \emph{Lie groupoid over} $\mathcal{F}$ is a groupoid over $\mathcal{F}$,
$F:\mathcal{E}\rightarrow\mathcal{F}$, such that the set of arrows
$F^{-1}(\mathit{Hom}(X,Y))$ is a manifold (for any objects $X$,
$Y$ of $\mathcal{F}$), the projection $F^{-1}(\mathit{Hom}(X,Y))\rightarrow\mathit{Hom}(X,Y)$
is a submersion and the composition of arrows in $\mathcal{E}$ is
smooth.

Now we can define the 2-categories $\mathcal{G}_{[n]}$ of generalized
Lie (super)groupoids. The objects of $\mathcal{G}_{[n]}$ are the
Lie supergroupoids over the supersemigroup $\left(\mathbb{R}^{0|n}\right)^{\mathbb{R}^{0|n}}$
(understood as a category with just one object; as we have seen above,
it is useful to identify it with the full subcategory of $\mathcal{S}$
with the object $\mathbb{R}^{0|n}$). Morphisms (generalized HS morphisms)
are Lie supergroupoids over $\left(\mathbb{R}^{0|n}\right)^{\mathbb{R}^{0|n}}\times\Delta_{1}$;
their compositions are defined as Lie supergroupoids over $\left(\mathbb{R}^{0|n}\right)^{\mathbb{R}^{0|n}}\times\Delta_{2}$.
The \emph{fibre} $\mathsf{G}_{0}$ of an object $\mathsf{G}$ of $\mathcal{G}_{[n]}$
is the fibre of $\mathsf{G}\rightarrow\left(\mathbb{R}^{0|n}\right)^{\mathbb{R}^{0|n}}$
over the unique object of $\left(\mathbb{R}^{0|n}\right)^{\mathbb{R}^{0|n}}$.

If $\mathbf{G}$ is a Lie supergroupoid on which
$\left(\mathbb{R}^{0|n}\right)^{\mathbb{R}^{0|n}}$ acts, it gives
us an object
$\left(\mathbb{R}^{0|n}\right)^{\mathbb{R}^{0|n}}\ltimes\mathbf{G}\rightarrow\left(\mathbb{R}^{0|n}\right)^{\mathbb{R}^{0|n}}$
of $\mathcal{G}_{[n]}$ (after all, objects of $\mathcal{G}_{[n]}$
are supergroupoids on which
$\left(\mathbb{R}^{0|n}\right)^{\mathbb{R}^{0|n}}$ acts up to
natural transformations); we will denote it for short
$\underline{\mathbf{G}}$ (notice that the fibre of
$\underline{\mathbf{G}}$ is $\mathbf{G}$). We get in this way an
embedding $\mathcal{G}\rightarrow\mathcal{G}_{[n]}$: for any
groupoid $\Gamma$, the corresponding object of $\mathcal{G}_{[n]}$
is $\underline{\Gamma^{\mathbb{R}^{0|n}}}$. Any object
$\mathsf{G}$ of $\mathcal{G}_{[n]}$ now defines a stack of all
1-morphisms
$\underline{M^{\mathbb{R}^{0|n}}}\rightarrow\mathsf{G}$; these
stacks will be called \emph{representable at the level} $n$.

We conclude with a concrete description of 1-morphisms
$\underline{M^{\mathbb{R}^{0|n}}}\rightarrow\mathsf{G}$. By
definition, we have to describe the superspace of arrows over
$\left(\mathbb{R}^{0|n}\right)^{\mathbb{R}^{0|n}}\times(0\rightarrow1)$
and their composition with arrows in
$\underline{M^{\mathbb{R}^{0|n}}}$ and in $\mathsf{G}$. However,
we can naturally identify this superspace of arrows with
$P\times\left(\mathbb{R}^{0|n}\right)^{\mathbb{R}^{0|n}}$, where
$P$ is the superspace of arrows over
$\mathit{id}\times(0\rightarrow1)$.

A morphism $\underline{M^{\mathbb{R}^{0|n}}}\rightarrow\mathsf{G}$
is thus equivalently given by a surjective submersion
$P\rightarrow M^{\mathbb{R}^{0|n}}$ and a right action of
$\mathsf{G}$ on $P$ satisfying two conditions: 1.~the map
$P\rightarrow M^{\mathbb{R}^{0|n}}$ is $\mathsf{G}$-equivariant
(where $\mathsf{G}$ acts on $M^{\mathbb{R}^{0|n}}$ via the
projection
$\mathsf{G}\rightarrow\left(\mathbb{R}^{0|n}\right)^{\mathbb{R}^{0|n}}$
and via the action of
$\left(\mathbb{R}^{0|n}\right)^{\mathbb{R}^{0|n}}$ on
$M^{\mathbb{R}^{0|n}}$) and 2.~the action of the fibre
$\mathsf{G}_{0}$ of $\mathsf{G}$ on $P$ makes $P\rightarrow
M^{\mathbb{R}^{0|n}}$ to a principal $\mathsf{G}_{0}$-bundle.

\subsection{Examples of stacks and of generalized Lie groupoids}

\subsubsection{{}``Categorified de Rham complex''}

This is one of the simplest possible examples. Let
$\mathbb{R}[k]\subset\Omega^{k}(\mathbb{R}^{0|1})$ denote the
(1-dimensional) group of closed $k$-forms on $\mathbb{R}^{0|1}$.
Since $\left(\mathbb{R}^{0|1}\right)^{\mathbb{R}^{0|1}}$ acts from
the right on $\mathbb{R}[k]$, it is an object of
$\mathcal{S}_{[1]}$, but since $\mathbb{R}[k]$ is a group, it
gives us an object of $\mathcal{G}_{[1]}$, namely
$\underline{\mathbb{R}[k]}$. An $\mathcal{S}_{[1]}$-morphism $\Pi
TM\rightarrow\mathbb{R}[k]$ is, tautologically, a closed $k$-form
on $M$. Here we shall consider 1-morphisms $\underline{\Pi
TM}\rightarrow\underline{\mathbb{R}[k]}$ in the 2-category
$\mathcal{G}_{[1]}$. As we explained above, to describe such a
morphism, we just have to describe the supermanifold $P$ of arrows
over $\mathit{id}\times(0\rightarrow1)$. The result is a principal
$\mathbb{R}[k]$-bundle $P\rightarrow\Pi TM$ in the category
$\mathcal{S}_{[1]}$,
i.e.~$\left(\mathbb{R}^{0|1}\right)^{\mathbb{R}^{0|1}}$ acts from
the right on $P$ and the maps $P\rightarrow\Pi TM$ and
$\mathbb{R}[k]\times P\rightarrow P$ are
$\left(\mathbb{R}^{0|1}\right)^{\mathbb{R}^{0|1}}$-equivariant.
Such $P$'s are easily seen to be classified by
$H^{k+1}(M,\mathbb{R})$ (see \cite{se2}).

Let us consider the exact sequence\[
0\rightarrow\mathbb{R}[k]\rightarrow\Omega^{k}(\mathbb{R}^{0|1})\rightarrow\mathbb{R}[k+1]\rightarrow0.\]
 All 1-morphisms $\underline{\Pi TM}\rightarrow\underline{\Omega^{k}(\mathbb{R}^{0|1})}$
are easily seen to be isomorphic. A lift of a 1-morphism
$\underline{\Pi
TM}\rightarrow\underline{\Omega^{k}(\mathbb{R}^{0|1})}$ to
$\underline{\Pi TM}\rightarrow\underline{\mathbb{R}[k]}$ is
equivalent to trivialization of the composed morphism
$\underline{\Pi
TM}\rightarrow\underline{\Omega^{k}(\mathbb{R}^{0|1})}\rightarrow\underline{\mathbb{R}[k+1]}$,
i.e. to an isomorphism (a 2-morphism) between the composition
$\underline{\Pi TM}\rightarrow\underline{\mathbb{R}[k+1]}$ and the
trivial morphism $\underline{\Pi
TM}\rightarrow\underline{\mathbb{R}[k+1]}$ (which is the same as a
trivialization of the corresponding principal
$\mathbb{R}[k+1]$-bundle over $\Pi TM$).

\subsubsection{Weil and Cartan models of equivariant cohomology }

We met this example in section
\ref{sub:Examples-of-approximations}, here we just rephrase it
from a different (and more natural) point of view. Let a Lie group
$G$ act on a manifold $N$; let us consider the stack of principal
$G$-bundles with equivariant maps to $N$, in other words, the
stack of HS morphisms $M\rightarrow\Gamma$, where $\Gamma$ is the
action groupoid $\Gamma$ of $G$ on $N$.

Objects of this stack have a lot of automorphisms. Let us consider
another, {}``rigidified'' stack of principal $G$-bundles with
equivariant maps to $N$ (as above) and with a choice of
connection. This stack is more rigid (any automorphism is
specified by its restriction to a point, if the base is
connected), but it is in the obvious sense homotopy equivalent to
the first one.

This rigidified stack is representable at level 1. To find the corresponding
object of $\mathcal{G}_{[1]}$, we have to take the groupoid whose
objects are triples $(P,h,\alpha)$, where $P\rightarrow\mathbb{R}^{0|1}$
is a principal $G$-bundle, $h:P\rightarrow N$ is an equivariant
map and $\alpha$ a connection on $P$. To get a small equivalent
groupoid, take $P$ to be the trivial bundle (any bundle over $\mathbb{R}^{0|1}$
can be trivialized): now the objects are pairs $(h,\alpha)$, $h$
is a map $\mathbb{R}^{0|1}\rightarrow N$ and $\alpha$ a $\frak{g}$-valued
1-form on $\mathbb{R}^{0|1}$, i.e.~they form the supermanifold $\Pi TN\times\Pi T\Pi\frak{g}$
(functions there are differential forms on $N$ with values in the
Weil algebra of $\frak{g}$). The invariant functions on $\Pi TN\times\Pi T\Pi\frak{g}$
are the basic subcomplex, i.e.~we got Weil model of equivariant cohomology.
Notice that we can get even smaller equivalent groupoid: we can always
make $\alpha$ to vanish at the origin by a gauge transformation.
By taking only these $\alpha$'s, we get Cartan model.

Equivariant cohomology of $N$ is the cohomology of the classifying
space of the action groupoid of $G$ on $N$. It seems probable that
one can get cohomology of classifying spaces of other (proper) Lie
groupoids by similar means, i.e.~by replacing the classifying
stack of HS morphisms $M\rightarrow\Gamma$ with a homotopy
equivalent rigid stack (whatever is the precise definition).
Unfortunately, we were not able to solve this problem, so we leave
it for later work.

\subsubsection{Quasi-Poisson groupoids}

As was noticed by A.~Vaintrob \cite{vain}, a Lie algebroid structure
on a vector bundle $A\rightarrow M$ is equivalent to an odd Poisson
structure on $\Pi A^{*}$ of degree 1. In this section we shall use
usual notation for $\mathbb{Z}$-graded supermanifolds, i.e.~we shall
denote $\Pi A^{*}$ as $A^{*}[1]$ (a $\mathbb{Z}$-grading can be
seen as an action of the {}``scaling'' 1-parameter group; degrees
are then weights with respect to this action). A Lie quasi-bialgebroid
structure on $A$ (as was observed in \cite{se1}) is equivalent to
a principal $\mathbb{R}[2]$-bundle $X\rightarrow A^{*}[2]$ in the
category of graded supermanifolds, with an $\mathbb{R}[2]$-invariant
odd Poisson structure on $X$; this Poisson structure can be projected
to $A^{*}[1]$, making $A$ to a Lie algebroid.

Lie quasi-bialgebroids are understood as infinitesimal objects corresponding
to quasi-Poisson groupoids. We can give a simple invariant definition
of the latter by integrating $X$ to an odd symplectic groupoid, keeping
track of the grading and of the $\mathbb{R}[2]$: a \emph{quasi-Poisson
structure} on a Lie groupoid $\Gamma\rightrightarrows M$ is a graded
odd symplectic groupoid $\breve{\Gamma}$ with symplectic form of
degree 1, with a submersion $\breve{\Gamma}\rightarrow\mathbb{R}[-1]$
that is a morphism of groupoids ($\breve{\Gamma}\rightarrow\mathbb{R}[-1]$
is the moment map for an action of $\mathbb{R}[2]$ on $\breve{\Gamma}$)
and with an isomorphism $\breve{\Gamma}/\!/\mathbb{R}[2]\cong T^{*}[1]\Gamma$.

Quasi-Poisson groupoids give us simple examples of objects of $\mathcal{G}_{[1]}$.
We only need a part of the structure of a quasi-Poisson groupoid: a graded Lie supergroupoid $\breve{\Gamma}$
over $\mathbb{R}[-1]$ such that the action of the scaling group on
the fibre of $\breve{\Gamma}$ can be extended to action of the scaling
semigroup $(\mathbb{R},\times)$. We get a groupoid over $\left(\mathbb{R}^{0|1}\right)^{\mathbb{R}^{0|1}}$
using the isomorphism $\left(\mathbb{R}^{0|1}\right)^{\mathbb{R}^{0|1}}=(\mathbb{R},\times)\ltimes\mathbb{R}[-1]$
(a general element of $\left(\mathbb{R}^{0|1}\right)^{\mathbb{R}^{0|1}}$
is of the form $\theta\mapsto a\theta+\beta$; $(\mathbb{R},\times)$
corresponds to $\theta\mapsto a\theta$ and $\mathbb{R}[-1]$ to $\theta\mapsto\theta+\beta$).

\lyxaddress{\noun{Dept.~of Theoretical Physics, Mlynsk\'a dolina F2, Comenius
University\\ 84248 Bratislava, Slovakia} }

{\texttt{kochan@sophia.dtp.fmph.uniba.sk \qquad severa@sophia.dtp.fmph.uniba.sk}}

\end{document}